\documentclass[12pt]{article}
\usepackage[dvips]{epsfig}
\usepackage{amsthm,amsmath,amssymb}
\usepackage{latexsym}
\usepackage{placeins}
\usepackage{color}
\usepackage[colorlinks]{hyperref}

\definecolor{blue}{rgb}{0,0,0.9}
\definecolor{red}{rgb}{0.9,0,0}
\definecolor{green}{rgb}{0,0.50,0.10}
\definecolor{violet}{rgb}{0.5804,0.0000,0.8275}

\newcommand{\green}[1]{\begin{color}{green}#1\end{color}}

\def\@themcountersep{}

\newtheorem{THEO}{Theorem}[section]
\newtheorem{ALGo}[THEO]{Algorithm}
\newtheorem{CONJ}[THEO]{Conjecture}
\newtheorem{COND}[THEO]{Condition}
\newtheorem{ASSUMP}[THEO]{Assumption}
\newtheorem{CORO}[THEO]{Corollary}
\newtheorem{DEFI}[THEO]{Definition}
\newtheorem{EXAMP}[THEO]{Example}
\newtheorem{FACT}[THEO]{Fact}
\newtheorem{HYPO}[THEO]{Hypothesis}
\newtheorem{LEMM}[THEO]{Lemma}
\newtheorem{PROB}[THEO]{Problem}
\newtheorem{PROP}[THEO]{Proposition}
\newtheorem{REMA}[THEO]{Remark}
\newcommand{\theo}{\begin{THEO}}
\newcommand{\algo}{\begin{ALGo} \rm}
\newcommand{\cond}{\begin{COND} \rm}
\newcommand{\assump}{\begin{ASSUMP} \rm}
\newcommand{\conj}{\begin{CONJ}}
\newcommand{\coro}{\begin{CORO}}
\newcommand{\defi}{\begin{DEFI} \rm}
\newcommand{\examp}{\begin{EXAMP} \rm}
\newcommand{\fact}{\begin{FACT}}
\newcommand{\hypo}{\begin{HYPO} \rm}
\newcommand{\lemm}{\begin{LEMM}}
\newcommand{\prob}{\begin{PROB} \rm}
\newcommand{\prop}{\begin{PROP}}
\newcommand{\rema}{\begin{REMA} \rm}
\newcommand{\etheo}{\end{THEO}}
\newcommand{\ealgo}{\end{ALGo}}
\newcommand{\econd}{\end{COND}}
\newcommand{\eassump}{\end{ASSUMP}}
\newcommand{\econj}{\end{CONJ}}
\newcommand{\ecoro}{\end{CORO}}
\newcommand{\edefi}{\end{DEFI}}
\newcommand{\eexamp}{\end{EXAMP}}
\newcommand{\efact}{\end{FACT}}
\newcommand{\ehypo}{\end{HYPO}}
\newcommand{\elemm}{\end{LEMM}}
\newcommand{\eprob}{\end{PROB}}
\newcommand{\eprop}{\end{PROP}}
\newcommand{\erema}{\end{REMA}}

\def\0{\mbox{\bf 0}}
\def\1{\mbox{\bf 1}}
\def\2{\mbox{\bf 2}}
\def\3{\mbox{\bf 3}}
\def\4{\mbox{\bf 4}}
\def\5{\mbox{\bf 5}}
\def\6{\mbox{\bf 6}}
\def\7{\mbox{\bf 7}}
\def\8{\mbox{\bf 8}}
\def\9{\mbox{\bf 9}}

\def\x{\mbox{\boldmath $x$}}
\def\y{\mbox{\boldmath $y$}}

\def\A{\mbox{\boldmath $A$}}

\def\H{\mbox{\boldmath $H$}}

\def\O{\mbox{\boldmath $O$}}

\def\Q{\mbox{\boldmath $Q$}}

\def\X{\mbox{\boldmath $X$}}
\def\Y{\mbox{\boldmath $Y$}}

\def\EC{\mbox{$\cal E$}}

\def\NC{\mbox{$\cal N$}}

\def\TC{\mbox{$\cal T$}}

\def\inprod#1#2{\langle#1, \, #2\rangle}

\def\Real{\mbox{$\mathbb{R}$}}
\def\coneK{\mbox{$\mathbb{K}$}}

\def\SymMat{\mbox{$\mathbb{S}$}}

\def\SymN{\mbox{$\mathbb{N}$}}

\def\DNN{\mbox{$\mathbb{DNN}$}}
\def\CPP{\mbox{$\mathbb{CPP}$}}

\def\bGamma{\mbox{$\bf{\Gamma}$}}

\def\inprod#1#2{\langle #1,\,#2\rangle}

\begin{document}

\title{ \Large 
%
%
%
Doubly nonnegative relaxations are equivalent to completely positive reformulations of 
quadratic optimization problems with block-clique graph structures
}

\bigskip
\author{
\normalsize 
Sunyoung Kim\thanks{Department of Mathematics, Ewha W. University, 52 Ewhayeodae-gil, Sudaemoon-gu, Seoul 03760, Korea 
			({\tt skim@ewha.ac.kr}). The research was supported
               by   NRF 2017-R1A2B2005119.}, \and \normalsize
Masakazu Kojima\thanks{Department of Industrial and Systems Engineering,
	Chuo University, Tokyo 192-0393, Japan 
({\tt kojima@is.titech.ac.jp}).
	This research was supported by Grant-in-Aid for Scientific Research (A) 26242027.},
 \and \normalsize
Kim-Chuan Toh\thanks{Department of Mathematics, and Institute of Operations Research and Analytics, National University of Singapore,
10 Lower Kent Ridge Road, Singapore 119076
({\tt mattohkc@nus.edu.sg}). 
This research is supported in part by the Ministry of Education, Singapore, Academic Research Fund (Grant number: R-146-000-257-112).
 } 
}
\date{\normalsize \today}

\maketitle 
\vspace*{-0.8cm}

\begin{abstract}
\noindent
We study the equivalence among a nonconvex QOP, its CPP and DNN relaxations 
under the assumption that the aggregated and correlative sparsity 
of the data matrices of the CPP relaxation is represented by a block-clique graph $G$. 
By exploiting the correlative sparsity, we  decompose the CPP relaxation problem into a clique-tree structured 
family of  smaller subproblems. 
Each subproblem is associated with a node of a clique tree of $G$. 
The optimal value can be obtained by applying an algorithm that we propose for solving the subproblems recursively 
from leaf nodes to the root node of the clique-tree. 
We establish the equivalence between the QOP and its DNN relaxation
from the equivalence between the reduced family of subproblems and their DNN relaxations
by applying the known results on: 
(i) CPP and DNN reformulation of a class of QOPs 
with linear equality, complementarity and binary constraints in $4$ nonnegative variables. (ii) DNN reformulation of 
a class of quadratically constrained convex QOPs with any size. (iii) DNN reformulation of LPs with any size. As a result, 
we show that 
a  QOP whose subproblems are the QOPs mentioned in (i), (ii) and (iii) 
 is equivalent to its DNN relaxation, 
if the subproblems form a clique-tree structured family induced from a block-clique graph. 
\end{abstract}


\vspace{0.5cm}


\section{Introduction}

Completely positive programming (CPP) relaxations of a class of quadratic optimization problems (QOPs) have received a great deal of attention as they can provide the  optimal value
of the original QOP 
\cite[etc.]{ARIMA2012,BOMZE2017,BURER2009,KIM2019}
even with their computational intractability. They are also referred as CPP reformulations of QOPs. 
As for computationally tractable alternatives, further relaxations of the CPP reformulations
to doubly nonnegative programming (DNN) relaxations
have been  studied in \cite{ARIMA2017,KIM2013}, and  
their effectiveness in obtaining good lower  bounds for the optimal value of 
the QOP have been demonstrated in  \cite{ITO2018,KIM2013}.
Since the equivalence among QOPs, their CPP  and DNN relaxations cannot be obtained in general, in this paper, we explore structured QOPs for which
the equivalence to their DNN reformulations can be established. In particular,
we focus on some structured sparsity  characterized by block-clique graphs \cite{JOHNSON1996} and 
partial convexity to establish their equivalence for a class of QOPs. 
As far as we are aware of, this is the first time that
the equivalence of nonconvex QOPs and their DNN relaxations are 
studied for this class of structured QOPs.

We start by introducing a conic optimization problem (COP) model to simultaneously represent 
QOPs, CPP problems and DNN problems.   
Let $\Real^n$ denote the $n$-dimensional Euclidean space and $\Real^n_+$ the
corresponding nonnegative orthant. 
We assume that each $\x  \in \Real^n$ is a column vector, and $\x^T$ denotes 
its transpose.
Let $\SymMat^n$ denote the space of $n \times n$ symmetric matrices.  
For every pair of $\A \in \SymMat^n$ and $\X \in \SymMat^n$, $\inprod{\A}{\X}$ stands for 
their inner product defined as the trace of $\A\X$. Let 
$I_{\rm eq}$ and $I_{\rm ineq}$ be disjoint finite subsets of positive integers and $\Q^p \in \SymMat^n$ $(p \in \{0\} \cup I_{\rm eq} \cup I_{\rm ineq})$. 
Given a closed (possibly nonconvex) cone $\coneK \subset  \SymMat^n$, 
we consider the following general COP:
\begin{eqnarray*}
\mbox{COP$(\coneK)$: } \zeta(\coneK) &=& \inf\left\{ \inprod{\Q^0}{\X} : 
\begin{array}{l}
\X \in \coneK, \ X_{11}=  1, \\ \inprod{\Q^p}{\X} = 0 \ (p \in I_{\rm eq}), \\ 
\inprod{\Q^p}{\X} \leq 0 \ (p \in I_{\rm ineq})
\end{array}
\right\} 
\end{eqnarray*}
The distinctive feature of COP$(\coneK)$ is that it only involves {\em homogeneous} equality and inequality constraints 
except $X_{11}=1$. 
A general equality standard form COP can be transformed to COP$(\coneK)$ with $I_{\rm ineq}=\emptyset$ 
in a straightforward manner (see Section 2.2). 
If $\bGamma^n = \{\x\x^T \in \SymMat^n : \x \in \Real^n_+ \}$ is chosen 
as the closed cone $\coneK \subset  \SymMat^n$, then
COP$(\bGamma^n)$ represents a QOP with quadratic equality and inequality constraints 
in $\x \in \Real^n_+$ (see Section 2.5). Notice that the inequality constraints 
$\inprod{\Q^p}{\X} \leq 0 \ (p \in I_{\rm ineq})$  
are dealt with separately from the equality constraints 
$\inprod{\Q^p}{\X} = 0 \ (p \in I_{\rm eq})$ without introducing slack variables. 
In particular, the inequality constraints 
$\inprod{\Q^p}{\X} \leq 0 \ (p \in I_{\rm ineq})$ with $\X\in \bGamma^n$ correspond to convex quadratic inequality constraints in Section 2.5. 
If $\coneK$ is the DNN cone of size $n$ (denoted as $\DNN^n$) 
or the CPP cone of size $n$ (denoted as $\CPP^n$), 
then COP$(\coneK)$ represents a general DNN or CPP  problem, respectively (see Section 2.2). 
They are known to serve as convex relaxations of the nonconvex QOP, COP$(\bGamma^n)$. In general, 
$\zeta(\DNN^n) \leq \zeta(\CPP^n) \leq \zeta(\bGamma^n)$ holds
since $\bGamma^n \subset \CPP^n \ \mbox{($=$ the convex hull of $\bGamma^n$)}  \subset \DNN^n$. 

The main purpose of this paper is to investigate 
structured sparsity, in particular, the aggregated sparsity \cite{FUKUDA2003} 
and the correlative sparsity \cite{KOBAYASHI2008} characterized by block-clique graphs \cite{JOHNSON1996},
 and 
partial convexity of the data matrices 
of COP$(\coneK)$ 
 to establish the equivalence among the three types of problems,  COP$(\bGamma^n)$, COP$(\CPP^n)$ and COP$(\DNN^n)$. 

Sparsity, especially the chordal graph sparsity, has been heavily used to improve the computational efficiency 
of solving semidefinite programming (SDP) problems.   
The sparsity exploitation technique 
\cite[etc.]{FUKUDA2003,KOBAYASHI2008,VANDENBERGE2014,WAKI2006} 
for SDP problems  
was based on the semidefinite matrix completion in order 
 to  reduce the size of the positive semidefinite variable matrix. More precisely,  it
replaces a large but sparse variable matrix  of a given SDP problem 
with smaller positive semidefinite variable matrices whose sizes are
determined by the maximal cliques of the extended chordal graph characterizing the aggregated sparsity of the data matrices of 
the SDP problem. On the other hand,  exploiting  sparsity in
CPP problems has not been studied in the 
literature, to the best of our knowledge, as  the studies on CPP problems have been mainly
for theoretical interests.   

We discuss the equivalence of COP$(\bGamma^n)$, COP$(\CPP^n)$ and COP$(\DNN^n)$  based on the following techniques and/or facts. \vspace{-2mm} 
\begin{description}
\item{(a) } The CPP and DNN matrix completion in \cite{DREW1998}. \vspace{-2mm} 
\item{(b) } Exploiting (aggregated and correlative) sparsity in chordal graph structured SDPs \cite{FUKUDA2003,KOBAYASHI2008}.\vspace{-2mm} 
\item{(c) } $\CPP^n = \DNN^n$ if $n \leq 4$.\vspace{-2mm} 
\item{(d)} CPP reformulation of a class of QOPs with linear equality, binary and complementarity constraints 
\cite[etc.]{ARIMA2012,BOMZE2017,BURER2009,KIM2019}.
\vspace{-2mm} 
\item{(e) } DNN reformulation of quadratically constrained convex QOPs in 
nonnegative variables. 
See Lemma~\ref{lemma:convexCOP} in Section 2.6.
\vspace{-2mm} 
\end{description}
We note that (c) is well-known.
A block-clique graph $G$ is a chordal graph in which any 
two maximal cliques intersect in at most one vertex \cite{JOHNSON1996}. 
It was shown in \cite{DREW1998} 
that every partial CPP (or DNN) matrix  whose specified entries are  determined by an undirected graph 
$G$  has a CPP (or DNN) completion if and only if $G$ is a block-clique graph.
 
Let $\coneK \in \{\bGamma^n,\CPP^n,\DNN^n\}$.  In our method,  the basic idea developed in (b) 
combined with (a), instead of positive semidefinite matrix completion \cite{GRON1984}, is applied to 
COP$(\coneK)$
with $\Q^p \in \SymMat^n$ $(p \in \{0\} \cup I_{\rm eq} \cup I_{\rm ineq})$.
The data structure of  COP$(\coneK)$ is characterized by 
a block-clique graph $G$ with the maximal cliques $C_r$ $(r=1,\ldots,\ell)$.
COP$(\coneK)$ is then decomposed into a family of 
$\ell$ smaller size subproblems 
according to a clique tree structure induced from $G$. 
The family of $\ell$ subproblems inherit the clique tree structure of the block-clique graph $G$. More precisely, 
 they are associated with the $\ell$ nodes of the clique tree.  Two distinct subproblems are almost 
 independent but weakly connected in the sense that they 
 share one scalar variable if they are adjacent in the clique tree and no common variable otherwise. In addition, 
each problem associated with a clique $C_r$ in the family is of the same form as COP$(\coneK)$, 
but the size of its matrix variable is decreased to the size of $C_r$. It is important to note that 
the decomposition is independent of
 the choice of $\coneK \in \{\bGamma^n,\CPP^n,\DNN^n\}$. 

We utilize the aforementioned decomposition  of COP($\coneK$) for two purposes:
to efficiently solve large scale COPs and to show the equivalence among COP$(\bGamma^n)$, COP$(\CPP^n)$ and 
COP$(\DNN^n)$.
For the first purpose of efficiently solving large scale COPs,  
 the optimal value $\zeta(\coneK)$  of 
COP($\coneK$) is computed by  solving the $\ell$ small decomposed subproblems in the family.
 We propose an algorithm for  sequentially
 solving the decomposed subproblems. As a result,  the proposed algorithm is more efficient
than directly solving the original COP($\coneK$) 
when its size  becomes increasingly large.
Here  we implicitly assume that $\coneK = \DNN^n$, although all the results would remain valid 
even when the decomposed subproblems are not numerically tractable.
We should emphasize that the decomposition into smaller subproblems is 
certainly beneficial computationally. For example, the decomposition of a DNN problem 
of size 1000  into $200$ DNN subproblems 
of size at most $10$ is certainly much more numerically tractable  than the original
DNN problem.

For the second purpose of showing the equivalence among COP$(\bGamma^n)$, COP$(\CPP^n)$ and 
COP$(\DNN^n)$, the decomposition is applied to a pair of COP$(\coneK_1)$ and COP$(\coneK_2)$ 
with two distinct $\coneK_1, \ \coneK_2 \in \{\bGamma^n,\CPP^n,\DNN^n\}$. 
Then, 
two families of subproblems, say family~$1$  from COP$(\coneK_1)$ and family~$2$ 
 from COP$(\coneK_2)$, are obtained. 
The equivalence between 
COP$(\coneK_1)$ and COP$(\coneK_2)$ is reduced to the equivalence of family~$1$ and family~$2$. 
We note that a pair of decomposed subproblems, 
one from family~$1$ and the other form family~2, 
associated with a common clique $C_r$ share the common objective function 
and constraints except for their cone constraints. 
Thus, (c), (d) and/or (e) can be applied for the equivalence of each pair. 
If all pairs of subproblems from families~$1$ and~$2$ are equivalent,
 then COP$(\coneK_1)$ and COP$(\coneK_2)$ are equivalent. 
 In particular, 
all of nonconvex QOPs with linear and complementarity 
constraints in $4$ variables, quadratically constrained convex QOPs with any size and LPs with 
any size can be included as subproblems in a single QOP
formulated as COP$(\bGamma)$,  which can then be 
equivalently reformulated as its DNN relaxation. 
We should mention that 
a  structured sparsity, 
{\it i.e.}, the aggregate and correlative sparsity characterized by a block-clique graph,  needs to be imposed on the QOP. 
Such a QOP may not appear  frequently in practice, but our study here has theoretical importance as nonconvex 
QOPs are NP hard to solve in general.

This paper is organized as follows: In Section 2, some basics on block-clique graphs,   
(a), (c), (d) and (e) are described. 
We also define the aggregate and correlative sparsity which are represented by an undirected graph. 
Sections 3 and 4 include the main results of this paper. 
In Section 3, the equivalence among COP$(\bGamma^n)$, COP$(\CPP^n)$ and 
COP$(\DNN^n)$ based on (c) and (d) 
is established by exploiting the aggregated sparsity characterized by block-clique graphs. 
In Section 4, we show how COP$(\coneK)$ with $\coneK \in \{\bGamma^n,\CPP^n,\DNN^n\}$ 
can be decomposed into smaller subproblems by exploiting the correlative sparsity.  
Then the equivalence results given in Section 3 as well as (e) are applied to a pair of subproblems 
induced from distinct COP$(\coneK_1)$ and COP$(\coneK_2)$ to  verify their equivalence 
as  mentioned above. 
We also present an algorithm to sequentially solve smaller decomposed subproblems. 
In Section 5, we illustrate two types of QOPs,  block-clique graph structured QOPs with linear equality and complementarity constraints and  partially convex QOPs, which can be formulated as 
the equivalent DNN problems. 
We conclude the paper in Section~6.

\section{Preliminaries}

\subsection{Notation and symbols}

We will consider the following cones in $\SymMat^n$.
\begin{eqnarray*}
\SymMat^n_+ &=& \mbox{the cone of $n \times n$ symmetric positive semidefinite matrices}, \\
\SymN^n &=& \mbox{the cone of $n \times n$ symmetric nonnegative matrices} \\
&=& \left\{ \X \in \SymMat^n : X_{ij} \geq 0 \ (i=1,\ldots,n, \ j=1,\ldots,n) \right\},  \\
\DNN^n & = & \mbox{ 
                         DNN cone}  \
           = \ \SymMat^n_+ \cap \SymN^n, \\ 
\bGamma^n & = & \left\{ \x\x^T \in \SymMat^n : \x \in \Real^n_+ \right\}, \\            
\CPP^n &=& \mbox{CPP  cone} \
= \ \mbox{the convex hull of }   \bGamma^n. 
\end{eqnarray*}
The following relation is well-known 
\cite{BERMAN2003}:
\begin{eqnarray}
\begin{array}{ll}
\CPP^{n} = \DNN^n 
& \mbox{if } n \leq 4, \\
\CPP^{n} \ \mbox{is a proper subset of } \DNN^n 
& \mbox{otherwise.}
\end{array}
\label{eq:inclusion}
\end{eqnarray} 

Let $N = \{1,\ldots,n\}$. 
For each nonempty subset $C$ of $N$, $\Real^C$ denotes the $|C|$-dimensional Euclidean space of column vectors of $x_i \ (i \in C)$, and 
$\SymMat^C$ the linear space of $|C| \times |C|$ 
symmetric matrices consisting of the elements 
$X_{ij}$ $(i,j) \in C \times C$. A vector in $\Real^C$ denoted by $\x_C$ is  regarded as a subvector of $\x \in \Real^n$, 
and a matrix in $\SymMat^C$ denoted by $\X_C$  as a principal submatrix of 
$\X \in \SymMat^n$. We define $\bGamma^C = \left\{\x_C\x_C^T : \x_C \in \Real^C_+\right\}$, and 
$\DNN^C$ and $\CPP^C$ to be the DNN cone and the CPP 
cone in $\SymMat^C$, respectively. 

We say a subset $F$ of $N \times N$ symmetric if it satisfies 
$(i,j) \in F \Leftrightarrow (j,i) \in F$.  For every symmetric subset $F$ of $N \times N$, we denote 
$\{(i,j) \in F : i\not=j\}$ and $F \cup \{(i,i): i \in N\}$ by $F^o$ and $\overline{F}$, respectively.

\subsection{Conversion of general COPs to COP$(\coneK)$}

Let $d_p \in \Real$ and 
$\overline{\Q}^p \in \SymMat^n$ $(p \in \{0\} \cup I_{\rm eq} \cup I_{\rm ineq})$. 
We assume $d_0=0$. 
For each $\coneK \in \{\bGamma^{n},\CPP^{n},\DNN^{n}\}$, 
consider a general COP in the standard equality form:  
\begin{eqnarray*}
\inf\left\{ \inprod{\overline{\Q}^0}{\Y} : 
\Y \in \coneK, \   \inprod{\overline{\Q}^p}{\Y} = d_p \ (p=1,\ldots,m)
\right\}.  
\end{eqnarray*}
Let
$ I_{\rm eq} = \{1,\ldots,m\}, \ I_{\rm ineq} = \emptyset$,  \  and 
$\Q^p = \begin{pmatrix} -d_p & \0^T \\ \0 & \overline{\Q}^p \end{pmatrix} \in \SymMat^{1+n} \ (p=0,1,\ldots,m).
$ 
Then $\inprod{\Q^p}{\X} = \inprod{\overline{\Q}^p}{\Y} - d_p X_{11}$ for every 
$\X = \begin{pmatrix} X_{11} & \x^T \\ \x & \Y \end{pmatrix} \in \SymMat^{1+n}$ $(p=0,1,\ldots,m)$. 
Therefore, the above standard equality form COP with $\coneK=\bGamma^n$, $\coneK=\CPP^n$ and $\coneK=\DNN^n$ 
can be rewritten as COP$(\bGamma^{1+n})$, COP$(\CPP^{1+n})$ and COP$(\DNN^{1+n})$, respectively. Conversely, if slack variables are introduced for the inequality constraints 
$\inprod{\Q^p}{\X} \leq 0$ $(p \in I_{\rm ineq})$ in COP$(\coneK)$, COP$(\coneK)$ 
can be converted into the standard equality form COP in a straightforward fashion.

\subsection{Chordal and block-clique graphs}

\begin{figure}
\begin{center}
 \ifpdf
 \includegraphics[width=0.38\textwidth]{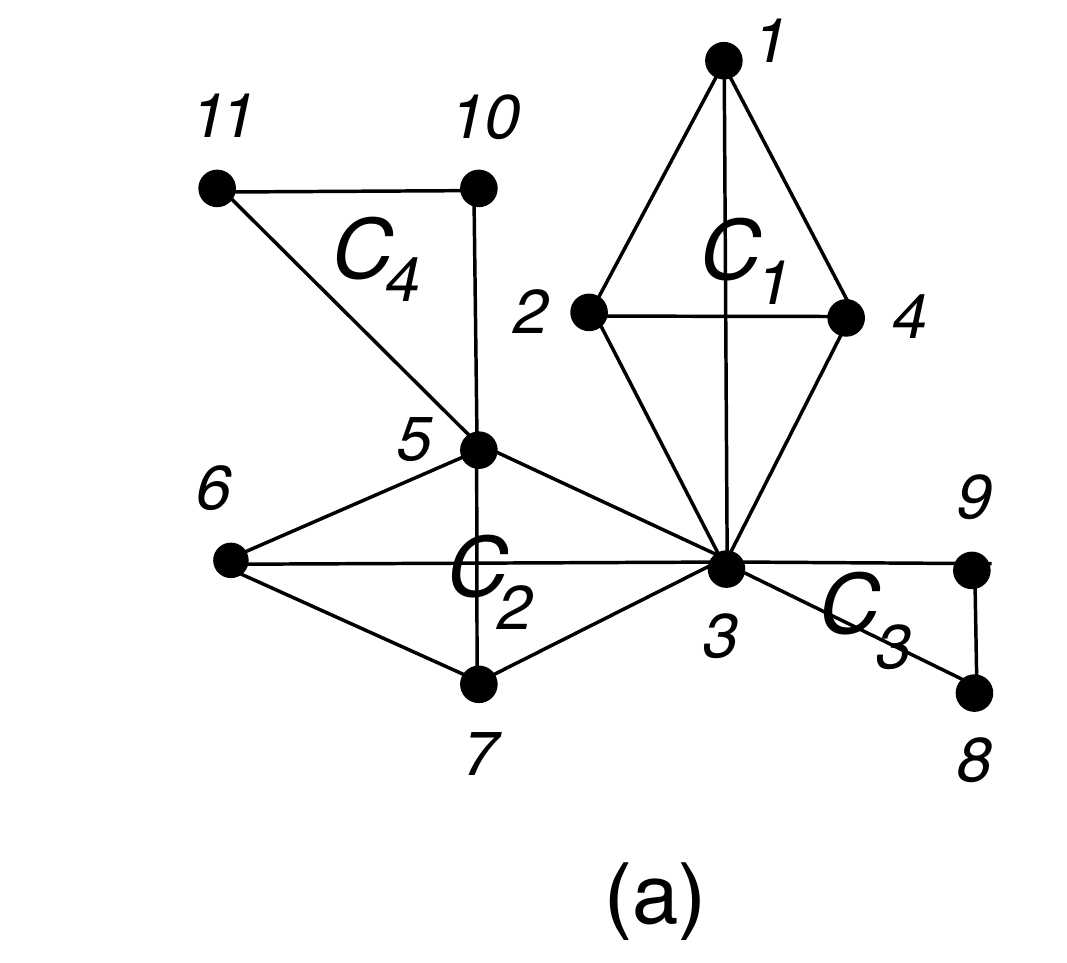}
 \else
 \includegraphics[width=0.38\textwidth]{blockClique_a.eps}
 \fi
 \ifpdf
\mbox{ \ }  \hspace{20mm} \mbox{ \ } 
 \includegraphics[width=0.36\textwidth]{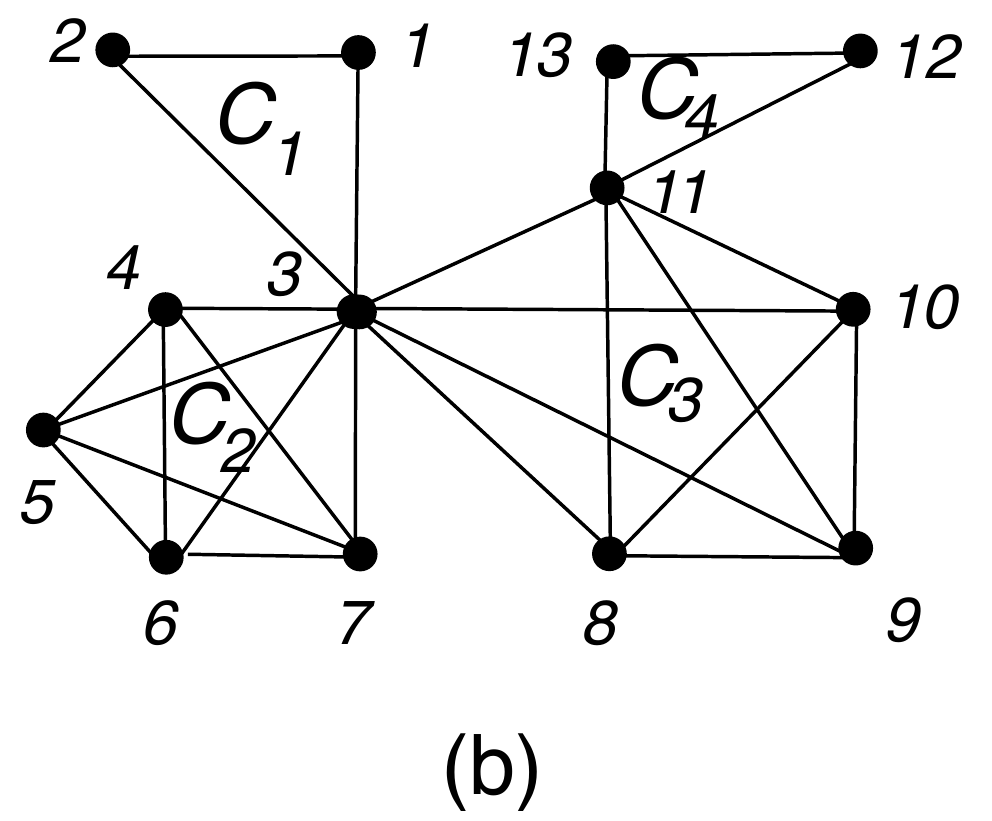}
 \else
 \includegraphics[width=0.36\textwidth]{blockClique_b.eps}
 \fi
\caption{
Illustration of block-clique graphs. 
In case (a), the maximal cliques are $C_1=\{1,2,3,4\}$, $C_2=\{3,5,6, 7\}$, $C_3=\{3,8,9\}$ and $C_4 = \{5,10,11\}$. 
In this case, every clique $C_q$ $(q=1,2,3,4)$ contains at most 4 nodes, 
so that the graph $G(N,E)$ satisfies the assumption of (ii) and (iii) of Theorem~\ref{theorem:main31}. 
In case (b), the maximal cliques are $C_1 = \{1,2,3\}$, $C_2=\{3,4,5, 6,7\}$, 
$C_3=\{3,8,9,10,11\}$ and $C_4=\{11,12,13\}$. If  a new edge $(10,12)$ is added, 
then  a new clique $C_5=\{10,11,12\}$ is created.
The resulting graph is no longer block-clique, but it remains to be chordal. 
}
\label{figure:Example1}
\end{center}
\end{figure}

We consider an undirected graph $G(N,E)$ with the node set $N = \{1,\ldots,n\}$ and the edge set $E$.  Here  
$E$ is a symmetric  subset of $\{ (i,j) \in N \times N : i \not= j \}$ (hence $E^o = E$) 
and $(i,j) \in E$ is identified with $(j,i) \in E$.
A graph $G(N,E)$ is called {\em chordal} if every cycle in $G(N,E)$ of length $4$ or more has 
a chord, and {\em block-clique} \cite{JOHNSON1996} if it is a chordal graph and any pair of two 
maximal cliques of $G(N,E)$ intersects in at most one node.  
See Figure 1 for examples of block-clique graphs.

Let $G(N,E)$ be a chordal graph with the maximal cliques $C_q$ $(q=1,\ldots,\ell)$. 
We assume that the graph is connected. 
If it is not,  
the subsequent discussion can be applied to each connected 
component. Consider an
undirected graph on the  maximal cliques, $G(\NC,\EC)$ with the node set $\{C_q : q=1,\ldots,\ell\}$ and the edge set 
$\EC = \{(C_{q},C_{r}) : C_{q} \cap C_{r} \not= \emptyset \}$. 
Since $G(N,E)$ is assumed to be connected, $G(\NC,\EC)$ is connected.
Then, add the weight $\left| C_{q} \cap C_{r} \right|$ to each 
edge $(C_{q},C_{r}) \in \EC$. Here $\left| C_{q} \cap C_{r} \right|$ denotes the number of 
nodes contained in the clique $C_{q} \cap C_{r} $. 
It is known that every maximum weight spanning tree $G(\NC,\TC)$ of $G(\NC,\EC)$ satisfies the 
following clique intersection property: 
\begin{eqnarray}
\left.
\begin{array}{l}
\mbox{for every pair of distinct cliques $C_{q}$ and $C_{r}$, $C_{q}\cap C_{r}$ is a subset of}\\
\mbox{every clique on the (unique) path connecting $C_{q}$ and $C_{r}$ in the tree.}
\end{array}
\right\} \label{eq:cip}
\end{eqnarray}
Such a tree is called as a {\it clique tree} of $G(N,E)$. We refer to \cite{BLAIR1993} for the fundamental 
properties of chordal graphs and clique trees including the clique intersection property and 
the running intersection property described below. 

\begin{figure}
\begin{center}
 \ifpdf
 \includegraphics[width=0.41\textwidth]{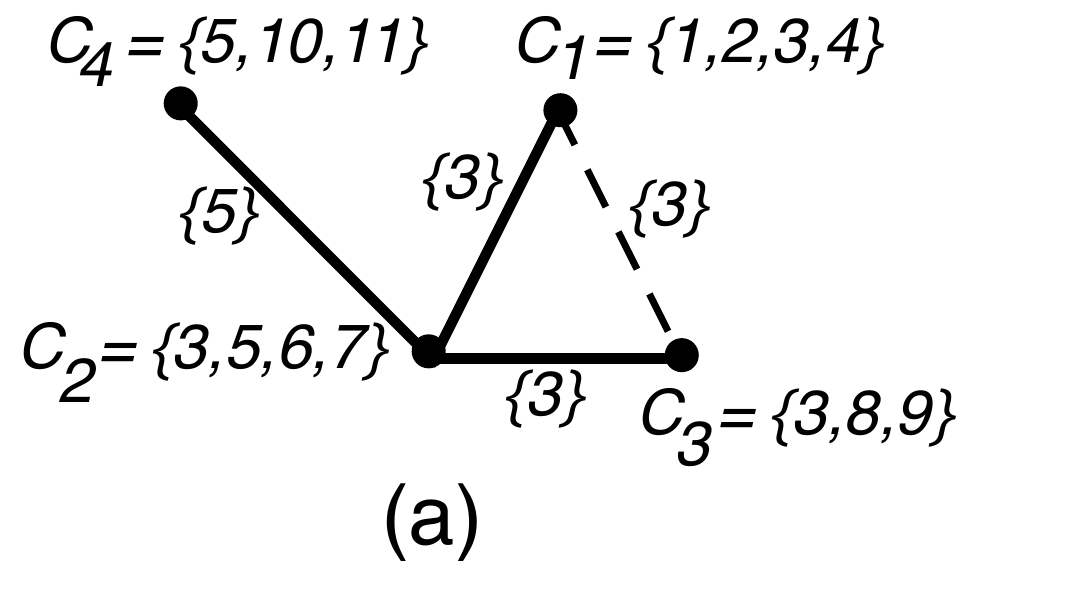}
 \else
 \includegraphics[width=0.41\textwidth]{CorrelativeSparsity_a.eps}
 \fi
 \ifpdf
\mbox{ \ }  \hspace{5mm} \mbox{ \ } 
 \includegraphics[width=0.48\textwidth]{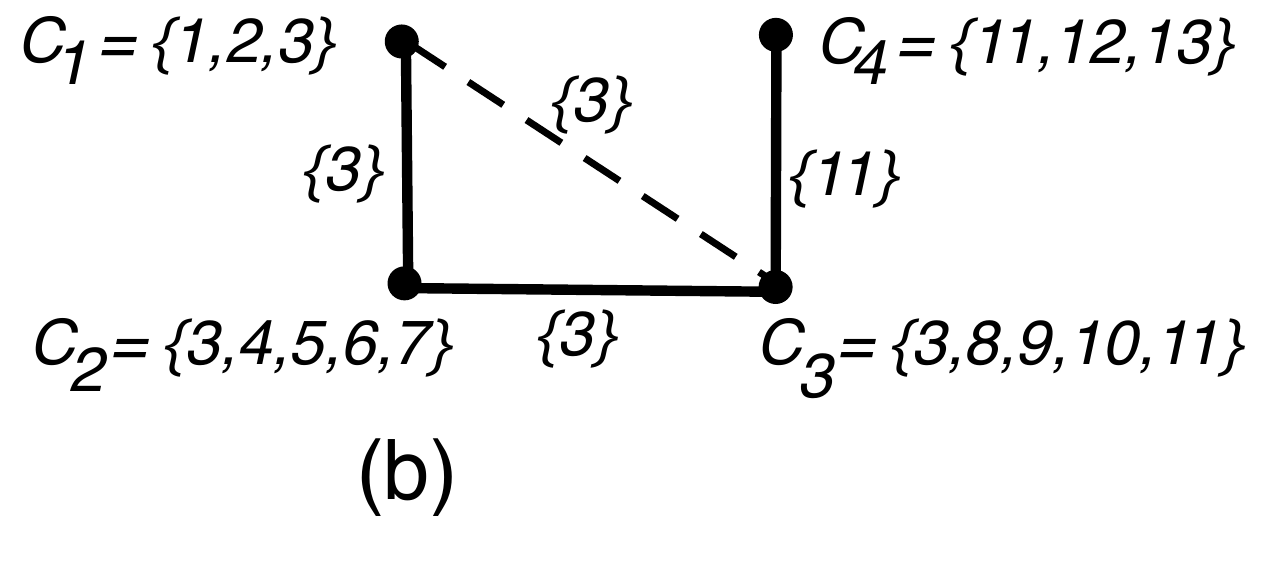}
 \else
 \includegraphics[width=0.48\textwidth]{CorrelativeSparsity_b.eps}
 \fi
\caption{
Illustration of clique trees. (a) and (b) are clique trees $G(\NC,\TC)$ of the block-clique graphs $G(N,E)$ 
shown in (a) and (b) of Figure 1, respectively. The notation $\{i\}$ added  on each edge $(C_q,C_r)$ denotes the intersection of the two cliques $C_q$ and $C_r$ of $G(N,E)$; 
for example, $(C_1,C_2)$ is an edge of the clique tree $G(\NC,\TC)$ in (a), and $C_1 \cap C_2 = \{3\}$. 
}
\label{figure:Example1}
\end{center}
\end{figure}

Now suppose that $G(N,E)$ is a connected block-clique graph. Then we know that 
$\left| C_{q} \cap C_{r} \right| = 1$ if $(C_{q},C_{r}) \in \EC$. Hence,  
every spanning tree of $G(\NC,\EC)$ is 
a clique tree. See Figure 2. Let $G(\NC,\TC)$ be a clique tree of $G(N,E)$. Choose an arbitrary maximal 
clique as a root node,  say $C_1$. 
The rest of the maximal cliques $C_2,\ldots,C_{\ell}$ 
can be renumbered such that for any pair of 
distinct $C_q$ and $C_r$, if $C_q$ is on the (unique) path from the root node $C_1$ to $C_r$ 
then $q < r$ holds by applying 
a topological sorting (ordering). In this way, the renumbered sequence of maximal cliques $C_1,C_2,\ldots,C_{\ell}$ satisfies {\em the running intersection property: }
\begin{eqnarray*}
& &
\mbox{$\forall r \in \{2,\ldots,\ell\}$, $\exists q \in \{1,\ldots,r-1\}$ such that }  
\left(C_{1} \cup \cdots \green{\cup} C_{r-1} \right) \cap C_r 
\subset C_{q}.  
\end{eqnarray*}
Let $r \in \{2,\ldots,\ell\}$, and $C_s$ the (unique) parent of $C_r$.   Then $s \in \{1,\ldots,r-1\}$ and 
\begin{eqnarray*}
\{k \} = C_s \cap C_r \subset \left(C_{1} \cup  \cdots \cup C_{r-1} \right) \cap C_r.  
\end{eqnarray*}
for some $k \in N$. If  $q \in  \{1,\ldots,r-1\}$ satisfies 
$ 
\left(C_{1} \cup \cdots \cup C_{r-1} \right)  \cap C_{r } \subset  C_{q}   
$ 
in the running intersection property above, then 
\begin{eqnarray*}
\{k \} = C_s \cap C_r \subset \left(C_{1} \cup \cdots \cup C_{r-1} \right) \cap C_r  \subset C_q  \cap  C_r = \{ k \}.  
\end{eqnarray*}
Here, the last equality follows from the assumption that $\left| C_q  \cap  C_r\right| \leq 1$. Therefore, we have shown the following 
result. 
\lemm \label{lemma:rip21}
 (The running intersection property applied to a block-clique graph.)
 Let $G(N,E)$ be a connected block-clique graph with the maximal cliques $C_1,C_2,\ldots,C_{\ell}$. 
Choose one of the maximal cliques arbitrary, say $C_1$. Then, the rest of the maximal cliques 
can be renumbered such that 
\begin{eqnarray}
\left. 
\begin{array}{l} 
\mbox{$\forall r \in \{2,\ldots,\ell\}$, $\exists q \in \{1,\ldots,r-1\}, \ \exists k_r \in C_r$ such that}  \\ 
\hspace{46mm} \left(C_{1}  \cup \cdots \cup C_{r-1} \right) \cap  C_r  = C_q \cap C_r = \{k_r\}  
\end{array}
\right\} \label{eq:rip21}
\end{eqnarray}
holds.
\elemm 

\subsection{Matrix completion}

We call an $n \times n$ matrix array $X_{ij}=X_{ji}$ $((i,j) \in N \times N)$ 
{\em a partial symmetric matrix} if a part of its elements $X_{ij}=X_{ji}$ $((i,j) \in F)$ are specified 
 for some symmetric $F \subset N \times N$ 
and the other elements are not specified. 
We denote a partial symmetric  matrix with specified elements 
$\bar{X}_{ij}=\bar{X}_{ji}$ $((i,j) \in F)$ by $[\bar{X}_{ij} : F]$. 
Given a property P characterizing 
a symmetric matrix in $\SymMat^n$ and a partial symmetric matrix $[\bar{X}_{ij} : F]$ 
for some symmetric $F \subset N \times N$, 
the matrix completion problem with property P is to find  values $\bar{X}_{ij} = \bar{X}_{ji}$  
 $((i,j) \not\in F)$ of unspecified elements $X_{ij} = X_{ji}$  
 $((i,j) \not\in F)$ such that the 
resulting $n \times n$ symmetric matrix $\bar{\X}$ has property P. 
We say that the partial symmetric matrix 
$[\bar{X}_{ij} :  F]$ has a {\em completion} $\bar{\X}$ with property P. 

We mainly consider CPP and DNN matrices in the subsequent discussions. 
In these matrices,  we may assume without loss of generality that 
$(i,i) \in F$ $(i \in N)$ 
since unspecified diagonal elements can be 
given as sufficiently large positive values to realize the property. 
Each partial symmetric  matrix $[\bar{X}_{ij} : F]$ can be associated with a graph $G(N,F^o)$. (Recall that $F^o=\{(i,j)\in F: i\not=j\}$.)
 Let $C_q$ 
$(q=1,\ldots,\ell)$ be the maximal cliques of $G(N,F^o)$. Then the partial symmetric  
matrix $[\bar{X}_{ij} : F]$ 
is decomposed into partial symmetric matrices $[\bar{X}_{ij} : C_q]$, which can be consistently  described 
as $\bar{\X}_{C_q}$,  $(q=1,\ldots,\ell)$. We say that a partial 
symmetric matrix $[\bar{X}_{ij} : F]$ is {\em partially 
DNN (partially CPP)} if every 
$\bar{\X}_{C_q}$ is DNN 
(CPP, respectively) in $\SymMat^{C_q}$ $(q=1,\ldots,\ell)$. 

\lemm \label{lemma:completion} \cite{DREW1998}
Let $F$ be a symmetric subset of $N \times N$ such that $(i,i) \in F$ for every $i \in N$. 
Then every partial CPP (DNN) matrix 
$[\bar{X}_{ij} :  F]$ has a CPP (DNN, respectively) 
completion, {\it i.e.}, a CPP (DNN, respectively) matrix $\X$ such that $X_{ij} = \bar{X}_{ij}$ $(i,j) \in F$, 
iff $G(N,F^o)$ is a block clique graph.
\elemm

\subsection{A class of QOPs and their CPP and DNN relaxations} 

Any quadratic function in nonnegative variables $x_2,\ldots,x_n$  can be represented as 
$\inprod{\Q}{\x\x^T}$ with $\x = (x_1,x_2,\ldots,x_n) \in \Real^n_+$ and $x_1^2 = 1$ for some $\Q \in \SymMat^n$. 
By introducing a variable matrix 
$\X \in \bGamma^n \equiv \{\x\x^T : \x \in \Real^n_+\}$,  the function can be rewritten 
as $\inprod{\Q}{\X}$ with $X_{11} = 1$. 
As a result, 
a general quadratically constrained QOP in nonnegative variables $x_2,\ldots,x_n$ can   also be represented 
as COP$(\bGamma^n)$ introduced in Section 1. Since $\bGamma^n \subset \CPP^n \subset \DNN^n$, 
$\zeta(\DNN^n) \leq \zeta(\CPP^n) \leq \zeta(\bGamma^n)$ holds in general. 

On the equivalence between QOPs and their CPP relaxations, Burer's
reformulation \cite{BURER2009} for a class of QOPs with linear constraints in nonnegative 
and binary variables is well-known. 
In this paper, we employ the following result, which is 
essentially equivalent to Burer's reformulation. See \cite{KIM2019} for CPP reformulations of more general class of QOPs.

\lemm \cite[Theorem 3.1]{ARIMA2018} \label{lemma:nonconvexCOP} For COP$(\bGamma^n)$,  assume that 
\begin{eqnarray*}
& & I_{\rm ineq} = \emptyset, \ 
\mbox{COP$(\bGamma^n)$ is feasible}, \\ 
& & \inprod{\Q^0}{\X} \geq 0 \ \mbox{if } \X \in \bGamma^n, \  X_{11} = 0 \ \mbox{and } \inprod{\Q^p}{\X} = 0 \ (p \in I_{\rm eq}), \\
& & \Q^p \in \SymMat^n_+ + \SymN^{n} \ \mbox{(the dual of $\DNN^n$)}  \ (p \in I_{\rm eq}). 
\end{eqnarray*}
Then, $\zeta(\bGamma^n) = \zeta(\CPP^n)$. 
\elemm

\examp \label{example:bbc} 
Let $\A$ be a $k \times n$ matrx, $I_{\rm comp} \subset \{(i,j) \in N \times N: 1 < i < j\}$ and $\Q^0 \in \SymMat^n$. Consider a QOP with linear equality and complementarity constraints in nonnegative 
variables $\x \in \Real^n_+$. 
\begin{eqnarray*}
\zeta_{\rm QOP} = \inf\left\{\x^T\Q^0\x : \x \in \Real^n_+, \ x_1 = 1, \ \A\x = \0, \ x_ix_j = 0 \ ((i,j) \in I_{\rm comp})\right\}.
\end{eqnarray*}
Define 
\begin{eqnarray*}
\Q^1 &=& \A^T\A \in \SymMat^n, \\
\Q^{ij} & = & \mbox{the $n\times n$ matrix with $1$ at the 
$(i,j)$th and $(j,i)$th elements,} \\
&& \mbox{and $0$ elsewhere.}
\end{eqnarray*}
Then, the QOP can be rewritten as 
\begin{eqnarray*}
\zeta_{\rm QOP} = \inf\left\{\inprod{\Q^0}{\X} : 
\begin{array}{l}
\X \in \bGamma^n, \ X_{11} = 1, \ \inprod{\Q^1}{\X} =0,  \\
\inprod{\Q^{ij}}{\X}= 0 \ ((i,j) \in I_{\rm comp})
\end{array}
\right\}.
\end{eqnarray*}
Enumerate $(i,j) \in I_{\rm comp}$ from $2$ through some integer
$m$, and let $I_{\rm eq} = \{1,\ldots,m\}$ and $I_{\rm ineq}=\emptyset$. Then QOP can be rewritten 
as COP$(\bGamma^n)$. Obviously $\Q^1 \in \SymMat^n_+$ and $\Q^{ij} \in \SymN^n$ $((i,j) \in I_{\rm comp})$. 
Hence, if the QOP is feasible and  
$ 
\O =\{\X \in \bGamma^m : X_{11} = 0, \ \inprod{\Q^1}{\X}= 0\} 
   =\{\x\x^T : \x\in\Real^n_+, \ x_1=0, \ \A\x = \0 \},
$ 
then all the assumptions in Lemma~\ref {lemma:nonconvexCOP} are satisfied. Consequently, 
$\zeta(\bGamma^n)=\zeta(\CPP^n)$ holds.
We note that the binary condition on variable $x$ can be represented as $x,y \geq 0$, $x+y=1$ and $xy=0$ with  a slack variable $y$. 
Thus, binary variables can be included in 
the QOP above. This QOP model covers various combinatorial optimization problems. See \cite{ARIMA2018,KIM2013} for more details.
\eexamp 

It is well-known that a convex QOP can be reformulated as 
an SDP (see, for example, \cite{FUJIE1997}). The following 
result may be regarded as a variation of the SDP reformulation to 
CPP and DNN reformulations of a convex QOP in nonnegative 
variables.

\lemm \label{lemma:convexCOP} 
Let $\widetilde{N} = N \backslash \{1\}$. 
Assume that $\Q^0_{\widetilde{N}} \in \SymMat^{n-1}_+$,  $\Q^p \in \SymMat^{n}_+$ $(p \in I_{\rm eq})$ and 
$\Q^p_{\widetilde{N}} \in \SymMat^{n-1}_+$ $(p \in I_{\rm ineq})$. 
Then $\zeta(\DNN^n) = \zeta(\CPP^n) = \zeta(\bGamma^n)$. Furthermore, if $\X = \begin{pmatrix} 1 & \y^T \\ \y & \Y \end{pmatrix} \in \DNN^n$ with some $\y \in \Real^{n-1}_+$ and $\Y \in \DNN^{n-1}$ is an optimal solution of 
{\rm COP}$(\DNN^n)$, then $\begin{pmatrix} 1 & \y^T \\ \y & \y\y^T \end{pmatrix}$ is a common optimal solution of 
{\rm COP}$(\bGamma^n)$, {\rm COP}$(\CPP^n)$ and {\rm COP}$(\DNN^n)$ with the objective value 
$\zeta(\DNN^n)= \zeta(\CPP^n)= \zeta(\bGamma^n)$. 
\elemm
\proof{The inequality $\zeta(\DNN^n) \leq \zeta(\CPP^n) \leq \zeta(\bGamma^n)$ 
follows from $\bGamma^n \subset \CPP^n \subseteq \DNN^n$. It suffices to show $\zeta(\bGamma^n) \leq \zeta(\DNN^n)$. 
Let  
$\X = \begin{pmatrix} 1 & \y^T \\ \y & \Y \end{pmatrix} \in \DNN^n$ with some 
$\y \in \Real^{n-1}_+$ and $\Y \in \DNN^{n-1}$ be an arbitrary feasible solution of 
COP$(\DNN^n)$. 
Then, 
\begin{eqnarray*}
& & \begin{pmatrix} 1 \\ \y \end{pmatrix} \in \Real^n_+, \ 
\overline{\X}  \equiv  \X - \begin{pmatrix} 0 & \0^T \\ \0 & \Y - \y\y^T \end{pmatrix}
= \begin{pmatrix} 1 & \y^T \\ \y & \y\y^T \end{pmatrix} 
\in \bGamma^n, \\ 
& & \inprod{\H^0}{\overline{\X}} =  1, \ 
\Q^p_{\widetilde{N}} \in \SymMat^{n-1}_+ \ (p \in \{0\}\cup  I_{\rm eq} \cup  I_{\rm ineq}), \ \Y - \y\y^T \in \SymMat^{n-1}_+, \\ 
& & \inprod{\Q^p}{\overline{\X}} = \inprod{\Q^p}{\X} - \inprod{\Q^p_{\widetilde{N}}}{\Y - \y\y^T} \leq \inprod{\Q^p}{\X}  \ (p \in \{0\}\cup  I_{\rm eq} \cup  I_{\rm ineq}). 
\end{eqnarray*}
From the last inequality, the assumptions and $\overline{\X} \in \CPP^n \subset \SymMat^n_+$, we see that 
\begin{eqnarray*}
& & 0 \leq \inprod{\Q^p}{\overline{\X}} \leq \inprod{\Q^p}{\X} = 0 \ \mbox{(hence } \inprod{\overline{\Q}^p}{\X} = 0 \mbox{)}   \ (p \in I_{\rm eq}), \\ 
& & \inprod{\Q^p}{\overline{\X}} \leq \inprod{\Q^p}{\X} \leq 0  \ (p \in  I_{\rm ineq}), \  
  \inprod{\Q^0}{\overline{\X}} \leq \inprod{\Q^0}{\X}. 
\end{eqnarray*}
Thus we have shown that $\overline{\X}$ is a feasible solution of COP$(\bGamma^n)$ whose 
objective value is not greater than that of the feasible solution $\X$ of COP$(\DNN^n)$. Therefore $\zeta(\bGamma^n) \leq \zeta(\DNN^n)$ has been shown. The second assertion follows by choosing an optimal solution of COP$(\DNN^n)$ for 
$\X$ in the  proof above. 
\qed
}

\examp \label{example:convexQOP}
(Quadratically constrained convex QOPs) 
Let $\widetilde{N}=\{2,\ldots,n\}$ and $I_{\rm ineq} = \{2,\ldots,m\}$. 
Let $\A$ be a $k \times n$ matrx and $\Q^p \in \SymMat^n$ be such that $\Q^0_{\widetilde{N}} \in \SymMat^{\widetilde{N}}_+$ 
$(p\in \{0\}\cup I_{\rm ineq})$. Consider a QOP:
\begin{eqnarray*}
\zeta_{\rm QOP} &=&\inf\left\{\x^T\Q^0\x : \x \in\Real^n_+, \ x_1= 1, \ \A\x=0, \ \x^T\Q^p\x \leq 0 \ (p \in I_{\rm ineq})\right\},
\end{eqnarray*}
which represents a general quadratically constrained convex QOP in  
nonnegative variables $x_i$ $(i=2,\ldots,n)$. If we let $\Q^1=\A^T\A \in \SymMat^n_+$ 
and $I_{\rm eq} = \{1\}$, we can represent the QOP as COP$(\bGamma^n)$. By Lemma~\ref{lemma:convexCOP}, not only 
$\zeta(\DNN^n) = \zeta(\CPP^n) = \zeta(\bGamma^n)$ holds, but also 
the DNN relaxation provides an optimal solution of the QOP.
\eexamp

\subsection{Two types of sparsity}

We consider two types of sparsity of the data matrices 
$\Q^p$ $(p \in \{0\} \cup I_{\rm eq} \cup I_{\rm ineq})$ 
of COP$(\coneK)$ with $\coneK \in \{\bGamma^n,\CPP^n,\DNN^n\}$, 
{\em the aggregated sparsity} \cite{FUKUDA2003} and 
{\em the correlative sparsity} \cite{KOBAYASHI2008}. 
We say that the aggregated sparsity of matrices $\A^p \in \SymMat^n$ $(p=1,\ldots,m)$ is 
represented by a graph $G(N,E)$ if 
$ 
\left\{(i,j) \in N \times N: A^p_{ij} \not= 0\right\} \subset  \overline{E} \equiv E \cup \{(i,i) : i \in N\} 
$ for every $p=1,\ldots,m$,  
and that their correlative sparsity is represented by a graph $G(N,E)$ with the maximal cliques $C_q$ $(q=1,\ldots,\ell)$ if 
\begin{eqnarray}
& & \mbox{$\forall p \in \{1,\ldots,m\}$, $\exists q \in \{q=1,\ldots,\ell\}$ such that } \nonumber \\ 
& &  \hspace{30mm} \{ (i,j) \in N \times N : A^p_{ij} \not= 0  \} 
\subset C_q \times C_q.
\end{eqnarray}
Note that if $ C_1,\ldots,C_{\ell}$ are the maximal cliques of $G(N,E)$,  then $\overline{E} = \cup_{q=1}^n C_q$.
We also note that a graph $G(N,E)$ which represents the aggregate (correlative) sparsity of 
$\A^p \in \SymMat^n$ $(p=1,\ldots,m)$ is not unique.

Our main interest in the subsequent discussion is 
a block-clique graph $G(N,E)$ which represents the aggregate (correlative) sparsity of $\A^p \in \SymMat^n$ $(p=1,\ldots,m)$. 
If we let 
$
F = \left\{(i,j) \in N \times N: A^p_{ij} \not= 0 \ \mbox{for some } p \right\} 
$, then $G(N,F^o)$ is ``the smallest'' graph which represents their aggregate sparsity. 
Since it is not block-clique in general,
a block-clique extension $G(N,E)$ of $G(N,F^o)$ is necessary. Assume that $G(N,F^o)$ is connected.
It is straightforward to verify that a node of  a connected block-clique graph is a cut node iff it is 
contained in at least two distinct maximal cliques of the graph. Therefore,
if there exists no cut node in $G(N,F^o)$, then the complete 
graph $G(N, (N \times N)^o)$ is the only block-clique extension of $G(N,F^o)$. Otherwise, take the maximal 
edge set $E \subset (N \times N)^o$ such that the 
graph $G(N,E)$ has the same cut nodes as $G(N,F^o)$. Then 
$G(N,E)$ forms the smallest block-clique extension of $G(N,F^o)$.

Obviously, if a graph $G(N,E)$ represents the correlative sparsity of matrices 
$\A^p \in \SymMat^n$ $(p=1,\ldots,m)$, then it also represents their aggregate sparsity. 
But the converse is not true as we see in the following example.
\examp
Let $n=3$, $N=\{1,2,3\}$ and 
\begin{eqnarray*}
\A^1=\begin{pmatrix}* & * & 0 \\ * & * & 0 \\ 0 & 0 & * \end{pmatrix}, \ 
\A^2= \begin{pmatrix} 0 & 0 & 0 \\ 0 & * & * \\ 0 & * & * \end{pmatrix}, \
\overline{\A}= \begin{pmatrix} * & * & 0 \\ * & * & * \\ 0 & * & * \end{pmatrix}
\end{eqnarray*}
where * denotes a nonzero element. We see that the aggregate sparsity pattern of the two matrices 
$\A^1$ and $\A^2$ corresponds to $\overline{\A}$. Thus their aggregated sparsity 
is represented by the graph with the maximal cliques $C_1 =\{1,2\}$ and $C_2=\{2,3\}$. 
But neither $C_1 \times C_1$ nor $C_2 \times C_2$ covers the nonzero elements of $\A^1$. 
We need to take the complete graph with the single maximal clique $N=\{1,2,3\}$ to represent 
the correlative sparsity of $\A^1$ and $\A^2$. 
\eexamp

\section{Exploiting the aggregated sparsity characterized by block-clique graphs}

Throughout this section, we assume that the aggregate sparsity of the data matrices 
$\Q^p$ $(p \in \{0\} \cup I_{\rm eq} \cup I_{\rm ineq})$ of COP($\coneK$) 
with $\coneK \in \{\bGamma^n,\CPP^n, \DNN^n\}$ 
is represented by a block-clique graph $G(N,E)$; 
\begin{eqnarray}
\{ (i,j) \in N \times N : Q^p_{ij} \not= 0 
\} \subset \overline{E}
\ (p \in \{0\} \cup I_{\rm eq} \cup I_{\rm ineq}).
\label{eq:sparsityPattern}
\end{eqnarray}
Under this assumption, we provide sufficient conditions for $\zeta(\DNN^n) = \zeta(\CPP^n)$ and 
$\zeta(\bGamma^n) = \zeta(\DNN^n) = \zeta(\CPP^n)$.

Note that 
the values of elements $X_{ij}$ $(i,j) \in \overline{E}$ determine
the value of $\inprod{\Q^p}{\X}$, {\it i.e.}, 
$\inprod{\Q^p}{\X} = \sum_{(i,j) \in \overline{E}}Q^p_{ij}X_{ij}$  
$(p \in \{0\} \cup I_{\rm eq} \cup I_{\rm ineq})$.  
All other elements $X_{ij}$ $(i,j) \not\in \overline{E}$ affect only 
 the cone constraint $\X \in \coneK$ in COP($\coneK$).  Thus,  the cone constraint can be replaced by 
 ``$[X_{ij}:\overline{E}]$ has a completion $\X \in \coneK$". 
More precisely, 
COP($\coneK$) can be written as 
\begin{eqnarray}
\zeta(\coneK) &=& \inf\left\{\sum_{(i,j) \in \overline{E}}Q^0_{ij}X_{ij} : 
\begin{array}{ll}
[X_{ij}: \overline{E}] \ \mbox{has a completion } \X \in \coneK, \ X_{11} = 1,\\
[3pt]
\displaystyle \sum_{(i,j) \in \overline{E}}Q^p_{ij}X_{ij} =0 \ (p \in I_{\rm eq}), \\ 
\displaystyle \sum_{(i,j) \in \overline{E}}Q^p_{ij}X_{ij} \leq 0 \ (p \in I_{\rm ineq})
\end{array}
\right\}.
\label{eq:COP32}
\end{eqnarray}

Let $C_q$ $(q=1,\ldots,\ell)$ be the 
maximal cliques of $G(N,E)$. We now introduce a (sparse) relaxation of COP \eqref{eq:COP32} ($=$ COP($\coneK$)). 
\begin{eqnarray}
\eta(\coneK) &=& \inf\left\{\sum_{(i,j) \in \overline{E}}Q^0_{ij}X_{ij} : 
\begin{array}{ll}
\X_{C_q}  \in \coneK^{C_q} \ (q=1,\ldots,\ell), \ X_{11} = 1, \\
[3pt]
\displaystyle \sum_{(i,j) \in \overline{E}}Q^p_{ij}X_{ij} =0 \ (p \in I_{\rm eq}), \\ 
\displaystyle \sum_{(i,j) \in \overline{E}}Q^p_{ij}X_{ij} \leq 0 \ (p \in I_{\rm ineq})
\end{array}
\right\}.
\label{eq:COP33}
\end{eqnarray}
Since we have been dealing with the case where 
$\coneK \in \{\bGamma^n, \CPP^n, \DNN^n\}$, 
$\coneK^{C_q}$ stands for $\bGamma^{C_q}$, $\DNN^{C_q}$ 
or $\CPP^{C_q}$ $(q=1,\ldots,\ell)$. In either  case, 
$\X \in \coneK$ implies $\X_{C_q} \in \coneK^{C_q}$ $(q=1,\ldots,\ell)$. 
Thus COP~\eqref{eq:COP33}  serves as a relaxation of 
COP~\eqref{eq:COP32}; $\eta(\coneK) \leq \zeta(\coneK)$.


\theo \label{theorem:main31} 
Let $\coneK \in \{\bGamma^n, \CPP^n, \DNN^n\}$. 
Assume that~\eqref{eq:sparsityPattern} holds for 
a block-clique graph $G(N,E)$ with the maximal cliques $C_q$ $(q=1,\ldots,\ell)$. Then, \vspace{-2mm} 
\begin{description}
\item{(i) } $\eta(\coneK) = \zeta(\coneK)$. \vspace{-2mm} 
\item{(ii) } If every clique $C_q$ contains at most $4$ nodes $(q=1,\ldots,\ell)$,  then 
COP~\eqref{eq:COP33} with $\coneK = \CPP^n$ and COP~\eqref{eq:COP33} with $\coneK = \DNN^n$ 
are equivalent; more precisely they share the same feasible solutions, the optimal solutions and the optimal value 
$\eta(\CPP^n) = \eta(\DNN^n).  
$\vspace{-2mm}
\item{(iii) }
If every clique $C_q$ contains at most $4$ nodes $(q=1,\ldots,\ell)$ and the assumptions of Lemma~\ref{lemma:nonconvexCOP} are satisfied, then $\eta(\bGamma^n) = \zeta(\bGamma^n) = \zeta(\CPP^n) = \eta(\CPP^n)=\eta(\DNN^n)$. \vspace{-2mm} 
\end{description}
\etheo
\proof{(i) Suppose that $\coneK \in \{\DNN^n,\CPP^n\}$. 
By Lemma~\ref{lemma:completion}, the condition 
``$[X_{ij}: (i,j)\in \overline{E}] \ \mbox{has a completion } \X \in \coneK$"
in COP~\eqref{eq:COP32} is equivalent to the condition 
``$\X_{C_q}  \in \coneK^{C_q} \ (q=1,\ldots,\ell)$" in COP~\eqref{eq:COP33}. Hence COPs~\eqref{eq:COP32} and~\eqref{eq:COP33} are equivalent, which implies that
$\eta(\coneK) = \zeta(\coneK)$. 
Now suppose that $\coneK = \bGamma^n$. Since we already know that $\eta(\bGamma^n) \leq \zeta(\bGamma^n)$, 
it suffices to show that $\zeta(\bGamma^n) \leq \eta(\bGamma^n)$. 
Choose a maximal clique that contains node $1$ 
among $C_1,\ldots,C_{\ell}$, say $C_1$. By Lemma~\ref{lemma:rip21},  we can renumber the rest of the maximal cliques 
such that the running intersection property~\eqref{eq:rip21} holds.  
Assume that $(\X_{C_1},\ldots,\X_{C_{\ell}})$ 
is a feasible solution of COP~\eqref{eq:COP33} with $\coneK=\bGamma^n$. We will construct an $\bar{\x} \in \Real^n_+$ such that 
$\overline{\X} \equiv\bar{\x}\bar{\x}^T$ is a feasible solution 
of COP$(\bGamma^n)$ with the same objective value 
of COP~\eqref{eq:COP33} with 
$\coneK=\bGamma^n$ at its feasible solution 
$(\X_{C_1},\ldots,\X_{C_{\ell}})$. 
Then $\zeta(\bGamma^n) \leq \eta(\bGamma^n)$ follows.
For each $r \in \{1,\ldots,\ell\}$, there exists an $\y_{C_r} \in \Real^C_+$ such that  $\X_{C_r} = \y_{C_r}\y_{C_r}^T$. 
For $r=1$, let $\bar{x}_i = [\y_{C_1}]_i$ $(i \in C_{1})$. By~\eqref{eq:rip21}, for each $r \in \{2,\ldots,\ell\}$,
 there exists a $k_r \in C_r$ such that $( C_1 \cup \cdots \cup C_{r-1})\cap C_r = \{k_r\}$. Hence, if we  define
 $\bar{x}_i =  [\y_{C_r}]_i \ (i \in C_r \backslash \{k_r\})$ for $r=2,\ldots,\ell$, then $\bar{\x} \in \Real^n_+$ 
 satisfies $\bar{x}_i =  [\y_{C_r}]_i \ (i \in C_r, r \in \{1,\ldots,\ell\})$ and $\overline{\X} = \bar{\x}\bar{\x}^T 
 \in \bGamma^n$ is a completion of $[X:\overline{E}]$. By construction, 
$\overline{\X}$  is a feasible solution of 
COP$(\bGamma^n)$, and attains the same objective value of COP~\eqref{eq:COP33} with 
$\coneK=\bGamma^n$ at its feasible solution 
$(\X_{C_1},\ldots,\X_{C_{\ell}})$.

(ii) By~\eqref{eq:inclusion}, 
$\CPP^{C_q}= \DNN^{C_q}$ $(q=1,\ldots,\ell)$, 
which implies the equivalence of COP~\eqref{eq:COP33} with $\coneK = \CPP^n$ and COP~\eqref{eq:COP33} 
with $\coneK = \DNN^n$.

(iii) The identity $\zeta(\bGamma^n)=\zeta(\CPP^n)$ follows from 
Lemma~\ref{lemma:nonconvexCOP}, and all other identities from 
(i) and (ii).}
\qed

\bigskip
To decompose the COP~\eqref{eq:COP33} based on the maximal cliques
$C_1,\ldots,C_\ell$, we first decompose
each $\Q^p$ $(p\in\{0\}\cup I_{\rm eq}\cup I_{\rm ineq})$ such that
\begin{eqnarray*}
\Q^p &=& \sum_{q=1}^{\ell} \widehat{\Q}^{pq}, \ 
\widehat{Q}^{pq}_{ij} = 0 \ \mbox{if } (i,j) \not\in C_q \times C_q \ (q=1,\ldots,\ell).
\end{eqnarray*}
(We note that the decomposition of $\Q^p$ above is not unique.)
Then COP~\eqref{eq:COP33} can be represented as
\begin{eqnarray}
\eta(\coneK) &=& \inf\left\{
\sum_{q=1}^{\ell} \inprod{\widehat{\Q}^{0q}_{C_q}}{\X_{C_q}} : 
\begin{array}{ll}
\X_{C_q}  \in \coneK^{C_q} \ (q=1,\ldots,\ell), \ X_{11} = 1, \\
[3pt]
\sum_{q=1}^{\ell} \inprod{\widehat{\Q}^{pq}_{C_q}}{\X_{C_q}} =0 \ (p \in I_{\rm eq}), \\
[3pt] 
\sum_{q=1}^{\ell} \inprod{\widehat{\Q}^{pq}_{C_q}}{\X_{C_q}} \leq 0 \ (p \in I_{\rm ineq})
\end{array}
\right\}. \label{eq:COPcliqueWise31}
\end{eqnarray}
\begin{REMA}
\begin{rm} \label{remark:LP}
(a) Note that
under the assumption~\eqref{eq:sparsityPattern}, if $\Q^p \in \SymMat^n_+$ 
(or $\Q^p \in \SymMat^n_+ + \SymN^n$), then  we can take $\widehat{\Q}^{pq} \in \SymMat^n_+$ 
(or $\widehat{\Q}^{pq} \in \SymMat^n_+ + \SymN^n$)  for the decomposition; see 
\cite[Theorem 2.3]{AGLER1988}. In this case, the equality constraints $\sum_{q=1}^{\ell} \inprod{\widehat{\Q}^{pq}_{C_q}}{\X_{C_q}} =0 \ (p \in I_{\rm eq})$ can be replaced by 
$\inprod{\widehat{\Q}^{pq}_{C_q}}{\X_{C_q}} =0 \ (q=1,\ldots,\ell,p \in I_{\rm eq})$ since 
$\inprod{\widehat{\Q}^{pq}_{C_q}}{\X_{C_q}} \geq 0$ for every 
$\X \in \coneK \in \{\bGamma^n,\CPP^n,\DNN^n\}$ is known. 
\\[5pt]
(b) Suppose that 
$C_q$ $(q=1,\ldots\ell)$ form a partition of $N$; $\cup_{q=1}^{\ell} C_q = N$ and $C_q \cap C_r = \emptyset$ $(q\not=r)$. 
Then we have 
$\widehat{\Q}^{pq}_{C_q} = \Q^{p}_{C_q}$ 
$(q=1,\ldots,\ell,p \in \{0\}\cup I_{\rm eq} \cup I_{\rm ineq})$.
In particular, if  $C_q = \{q\}$ $(q=1,\ldots,\ell=n)$, then 
COP~\eqref{eq:COPcliqueWise31} turns out to be an LP  of the form 
\begin{eqnarray*}
\eta(\coneK) = \inf \left\{
\sum_{q=1}^{\ell} Q^0_{ii}X_{ii} : 
\begin{array}{l}
X_{ii} \in \Real_+, \ (i=1,\ldots,n), X_{11} = 1, \\
\sum_{i=1}^{n} Q^p_{ii}X_{ii} = 0 \ (p \in I_{\rm eq}), \\
\sum_{i=1}^{n} Q^p_{ii}X_{ii} \leq 0 \ (p \in I_{\rm ineq})
\end{array}
\right\}. 
\end{eqnarray*}
Here $\bGamma^1 = \CPP^1 = \DNN^1 =\Real_+$. Although  
this observation itself is trivial, it is important in the sense that 
LPs can be embedded as subproblems in a QOP that can be reformulated 
as its DNN relaxations, as we will see in (v) of Theorem~\ref{theorem:main45} and Section 5.2.
\end{rm}
\end{REMA}

\section{Exploiting correlative sparsity characterized by block-clique graphs}

We consider COP$(\coneK)$ with $\coneK \in \{\bGamma^n, \CPP^n, \DNN^n\}$. 
Throughout this section,  we assume that  
the correlative sparsity \cite{KOBAYASHI2008} of the data matrices $\Q^p$  
$(p \in  I_{\rm eq} \cup I_{\rm ineq})$ of the linear equality and inequality constraints of 
COP$(\coneK)$ 
is represented by a connected block-clique graph $G(N,E)$ with the maximal cliques $C_q$ $(q=1,\ldots,\ell)$
and that the aggregate sparsity of $\Q^0$ is represented by the same block-clique graph $G(N,E)$, i.e.,
\begin{eqnarray}
& & \mbox{$\forall p \in I_{\rm eq} \cup I_{\rm ineq}$, $\exists q \in \{q=1,\ldots,\ell\}$ such that } \nonumber \\ 
& &  \hspace{40mm} \{ (i,j) \in N \times N : Q^p_{ij} \not= 0  \} 
\subset C_q \times C_q, \label{eq:corrSparsity31}\\
& & \{ (i,j) \in N \times N : Q^0_{ij} \not= 0  \} \subset \overline{E} \equiv \cup_{i=1}^{\ell} C_q \times C_q.    \label{eq:corrSparsity32}
\end{eqnarray}
In addition, we assume that 
COP$(\coneK)$ has an optimal solution. 
It should be noted that~\eqref{eq:corrSparsity31} and~\eqref{eq:corrSparsity32} 
imply~\eqref{eq:sparsityPattern}, thus, all the results discussed in the previous section remain valid. 
In particular, by (i) of Theorem~\ref{theorem:main31}, COP$(\coneK)$ is equivalent to 
COP~\eqref{eq:COP33}. 

Choose a clique that contains node $1$ among $C_q$ $(q=1,\ldots,\ell)$, say $C_1$. 
By Lemma~\ref{lemma:rip21}, 
the maximal cliques $C_2,\ldots,C_{\ell}$ can be renumbered 
such that the running intersection property~\eqref{eq:rip21} holds. 
Let $F_r = \cup_{q=1}^r C_q \times C_q$ $(r=1,\ldots,\ell)$. 
We then obtain from~\eqref{eq:rip21} that  
\begin{eqnarray}
& & \forall r \in \{2,\ldots,\ell\}, \exists k_r \in C_r \; \mbox{such that}\nonumber \\
& & \hspace{20mm}  F_{r-1}\cap (C_r \times C_r) = \left( \cup_{q=1}^{r-1} C_{q} \times C_{q} \right) \cap (C_r \times C_r) = \{ (k_r,k_r) \}. 
\label{eq:rip31} 
\end{eqnarray}
By~\eqref{eq:corrSparsity31}, there exists a partition $L_1,\ldots,L_{\ell}$ of 
$I_{\rm  eq} \cup I_{\rm  ineq}$
({\it i.e.}, $\cup_{q=1}^{\ell} L_q = I_{\rm  eq} \cup I_{\rm  ineq}$ and $L_{q} \cap L_{r} = \emptyset$ if $1 \leq q < r \leq \ell$) such that  
$ 
\{ (i,j) \in N \times N : Q^p_{ij} \not= 0  \}  \subset C_q \times C_q \ \mbox{if $p \in L_q$ } \ (q=1,\ldots,\ell). 
$ 
Hence, we can rewrite COP~\eqref{eq:COP33}, which is equivalent to COP$(\coneK)$,  as 
\begin{eqnarray}
\mbox{$\overline{\rm COP}(\coneK)$: } 
\bar{\zeta}(\coneK) &=& \inf\left\{\sum_{(i,j) \in F_{\ell}}Q^0_{ij}X_{ij} : 
\begin{array}{ll}
\X_{C_q}  \in \coneK^{C_q} \ (q=1,\ldots,\ell), \ X_{11} = 1,\\
[3pt]
 \displaystyle \sum_{(i,j) \in C_q \times C_q} Q^p_{ij}X_{ij} = 0 \\ 
 \hspace{15mm} (p \in I_{\rm eq}^q, \ q=1,\ldots,\ell) ,\\ 
 \displaystyle \sum_{(i,j) \in C_q \times C_q} Q^p_{ij}X_{ij} \leq 0 \\ 
\hspace{15mm}   (p \in  I_{\rm ineq}^q, \ q=1,\ldots,\ell)
\end{array}
\right\}. \quad
\label{eq:COPbar}
\end{eqnarray}
Here $I_{\rm eq}^q = L_q\cap I_{\rm eq}\ \mbox{and } \ I_{\rm ineq}^q  = \ L_q\cap I_{\rm ineq} 
\ (q=1,\ldots,\ell)$. 
In particular, $(\overline{\X}_{C_1},\ldots,\overline{\X}_{C_{\ell}})$
is an optimal solution of $\overline{\rm COP}(\coneK)$ with the objective 
value $\bar{\zeta}(\coneK) = \zeta(\coneK)  > -\infty$ if 
$\overline{\X}$ is an optimal solution of COP$(\coneK)$.  

\subsection{Recursive reduction of 
\mbox{$\overline{\rm COP}(\coneK)$} to a sequence of smaller-sized COPs}

We first describe the basic idea on
how to reduce $\overline{\rm COP}(\coneK)$ to smaller COPs 
by recursively eliminating variable matrices $\X_{C_r}$ $(r=\ell,\ldots,2)$. 
The constraints of $\overline{\rm COP}(\coneK)$ except $X_{11}=1$
can be decomposed into $\ell$ families of constraints
\begin{eqnarray}
\sum_{(i,j) \in C_q \times C_q} Q^p_{ij}X_{ij} = 0 \
 (p \in I_{\rm eq}^q) \ \mbox{and }
\sum_{(i,j) \in C_q \times C_q} Q^p_{ij}X_{ij} \leq 0 \
(p \in  I_{\rm ineq}^q) \label{eq:decompConst}
\end{eqnarray}
in the variable matrix $\X_{C_q}  \in \coneK^{C_q}$ $(q=\ell,\ldots,1)$. Although they may look  independent, they are 
weakly connected in the sense that
\begin{eqnarray*}
(\X_{C_1},\ldots,\X_{C_{r-1}})\ \mbox{and } \X_{C_r} \ 
\mbox{share only one scalar variable } X_{k_rk_r} 
\end{eqnarray*}
$(r=\ell,\ldots,2)$. This relation follows from~\eqref{eq:rip31}. 

Let $r=\ell$. Then,
the linear objective function 
of $\overline{\rm COP}(\coneK)$ can be decomposed into two 
non-interactive terms such that
\begin{eqnarray}
\sum_{(i,j) \in F_{\ell-1}} Q^0_{ij}X_{ij} + \sum_{(i,j)\in C_{\ell} \times C_{\ell} \backslash \{(k_{\ell},k_{\ell})\}} Q^0_{ij}X_{ij}.
\label{eq:decompObj}
\end{eqnarray}
Hence, if  the value of $X_{k_{\ell}k_{\ell}} \geq 0$ is specified,
the subproblem of minimizing the second term $\sum_{(i,j)\in C_{\ell} \times C_{\ell} \backslash \{(k_{\ell},k_{\ell})\}} Q^0_{ij}X_{ij}$ over the decomposed 
constraint~\eqref{eq:decompConst} with $q=\ell$ in 
$\X_{C_{\ell}} \in \coneK^{C_{\ell}}$  can be solved 
independently from $\overline{\rm COP}(\coneK)$. 
(This subproblem  corresponds to 
$\widetilde{\mbox{\rm P}}_{\ell}(\coneK^{C_{\ell}},X_{k_rk_r})$ 
defined later in \eqref{eq:Prr}). 
Since all equalities and inequalities in the constraint~\eqref{eq:decompConst} are homogeneous, the 
optimal value is proportional to the specified value 
$X_{k_{\ell}k_{\ell}} \geq 0$, i.e., equal to 
$\tilde{\eta}_{\ell}(C_{\ell},1)X_{k_{\ell}k_{\ell}}$, 
where $\tilde{\eta}_{\ell}(C_{\ell},1)$ denotes 
the optimal value of the subproblem with $X_{k_{\ell}k_{\ell}}$ 
specified to $1$. Thus, 
 the constraint~\eqref{eq:decompConst} with $q=\ell$ in 
the varible matrix $\X_{C_{\ell}}$ can be eliminated from 
$\overline{\rm COP}(\coneK)$ by replacing the objective 
function~\eqref{eq:decompObj} with 
\begin{eqnarray*}
\sum_{(i,j) \in F_{\ell-1}} Q^0_{ij}X_{ij} 
+ \tilde{\eta}_{\ell}(C_{\ell},1) X_{k_rk_r} = \sum_{(i,j) \in F_{\ell-1}} Q^{0\ell-1}_{ij}X_{ij}, 
\end{eqnarray*}
where $Q^{0\ell-1}_{ij}= Q^0_{ij}
+ \tilde{\eta}_{\ell}(C_{\ell},1)$ 
if $i=j=k_r$ and $Q^{0\ell-1}_{ij}= Q^0_{ij}$ otherwise. 
(The resulting problem will be represented as 
P$_{\ell-1}(\coneK)$ in later discussion). 
We continue this elimination process and
updating $Q^{0r}_{ij} $ 
to $ Q^{0r-1}_{ij}$ for $r=\ell-1,\ldots2$, till we obtain 
\begin{eqnarray*}
\inf\left\{\sum_{(i,j) \in C_1 \times C_1} Q^{01}_{ij}X_{ij}: 
\begin{array}{l}
\X_{C_1} \in \coneK^{C_1}, \ X_{11} = 1, \\[3pt]
\displaystyle \sum_{(i,j) \in C_1 \times C_1}Q^{p}_{ij}X_{ij} = 0 \ (p\in I^1_{\rm eq}), \\[3pt]
\displaystyle \sum_{(i,j) \in C_1 \times C_1}Q^{p}_{ij}X_{ij} \leq 0 \ (p\in I^1_{\rm ineq})
\end{array}
\right\}, 
\end{eqnarray*}
which has the same optimal value as $\overline{\rm COP}(\coneK)$. 
(The above problem corresponds to P$_1(\coneK)$ 
to be defined later). 

In the above brief description of the elimination process, we have implicitly 
assumed that the subproblem with $X_{k_rk_r}$ specified to $1$ has
an optimal solution, but it may be infeasible or unbounded. 
We need to deal with such cases. In the subsequent discussion,  we also show how an 
optimal solution of $\overline{\rm COP}(\coneK)$ is retrieved in detail.

To embed a recursive structure in $\overline{\rm COP}(\coneK), $
we introduce some notation.
Let $k_1 = 1$, and let 
\begin{eqnarray*}
\Phi_r(\coneK) & = & \left\{(\X_{C_1},\ldots,\X_{C_r}) : \begin{array}{ll}
\X_{C_q}  \in \coneK^{C_q} \ (q=1,\ldots,r), \ X_{11} = 1,\\
 \displaystyle \sum_{(i,j) \in C_q \times C_q} Q^p_{ij}X_{ij} = 0 \
(p \in I_{\rm eq}^q, q=1,\ldots,r) \\ 
 \displaystyle \sum_{(i,j) \in C_q \times C_q} Q^p_{ij}X_{ij} \leq 0 \
(p \in I_{\rm ineq}^q, q=1,\ldots,r)
\end{array}
\right\}, \\[5pt]
\Psi_r(\coneK^{C_r},\lambda) & = & \left\{ \X_{C_r}\in \coneK^{C_r} : 
\begin{array}{l}
 X_{k_rk_r} = \lambda,\\
  \displaystyle \sum_{(i,j) \in C_r \times C_r} Q^p_{ij}X_{ij} = 0 \ (p \in I_{\rm eq}^r)\\ 
  \displaystyle \sum_{(i,j) \in C_r \times C_r} Q^p_{ij}X_{ij} \leq 0 \ (p \in I_{\rm ineq}^r)\\ 
\end{array}
\right\}, 
\end{eqnarray*}
$(r=\ell,\ldots,1)$. 
By~\eqref{eq:rip31}, each $\Phi_r(\coneK)$  $(r=\ell,\ldots,2)$ can be represented as 
\begin{eqnarray}
\Phi_r(\coneK) & = & \left\{ (\X_{C_1},\ldots,\X_{C_r}) : 
\begin{array}{l}  (\X_{C_1},\ldots,\X_{C_{r-1}}) \in \Phi_{r-1}(\coneK), \\ 
\X_{C_r} \in \Psi_r(\coneK^{C_r},X_{k_rk_r}) 
\end{array}
\right\}.  \label{eq:recursive} 
\end{eqnarray}
This recursive representation of $\Phi_r(\coneK)$  $(r=\ell,\ldots,2)$ plays an essential role in the 
discussions below. 

We now construct a sequence of COPs: 
\begin{eqnarray}
\mbox{P$_{r}(\coneK)$:} \quad \eta_r(\coneK) = \inf \left\{ \sum_{(i,j) \in F_r} Q^{0r}_{ij}X_{ij} : 
(\X_{C_1},\ldots,\X_{C_r}) \in \Phi_r(\coneK) 
\right\}
\label{eq:Pr}
\end{eqnarray}
$(r=\ell,\ldots,1)$ 
such that 
\begin{eqnarray}
\left. 
\begin{array}{l}
\mbox{if $\overline{\X}$ is an optimal solution of COP$(\coneK)$, then 
$(\overline{\X}_{C_1},\ldots,\overline{\X}_{C_{q}})$ is} \\
\mbox{an optimal solution of P$_{q}(\coneK)$ with the optimal value $\eta_{q}(\coneK)=\bar{\zeta}(\coneK) > -\infty$}
\end{array}
\right\}, \label{eq:equivalence31}
\end{eqnarray}
where $Q^{0r}_{ij} \in \Real \cup \{ \infty \}$  $((i,j) \in F_q)$ $(r=\ell,\ldots,1)$. 
We note that  every $Q^{0q}_{ij}$ is fixed to $Q^{0}_{ij} \in \Real$ 
for $i\not=j$, but $Q^{0q}_{ii} \in \Real \cup \{ \infty \}$ $((i,i) \in F_q)$ $(q=\ell,\ldots,1)$ 
are updated in the sequence. If we assign $Q^{0q}_{ii} = +\infty$ for some $(i,i) \in F_q$, then the objective 
quadratic function $\sum_{(i,i) \in F_q} Q^{0q}_{ij}X_{ij}$ takes $\infty$ at $(\X_{C_1},\ldots,\X_{C_q}) \in \Phi_q(\coneK)$ 
unless $X_{ii} = 0$. Thus, if  $Q^{0q}_{ii} = +\infty$, then we must  take $X_{ii} = 0$ in P$_q(\coneK)$. 

As $r$ decreases from $\ell$ to $1$ in the sequence,  
we obtain at the final iteration:
\begin{eqnarray*}
\mbox{P$_{1}(\coneK)$:  } \eta_1(\coneK) 
& = & \inf \left\{ \sum_{(i,j) \in C_1 \times C_1} Q^{0r}_{ij}X_{ij} :  \X_{C_1} \in \Phi_1(\coneK) \right\},  
\end{eqnarray*}
which is equivalent to $\overline{\rm COP}(\coneK)$, {\it i.e.}, 
\eqref{eq:equivalence31} holds for $q = 1$.

To  construct the sequence P$_r(\coneK)$ $(r=\ell,\ldots,1)$, 
we first set 
\begin{eqnarray}
Q^{0\ell}_{ij} = Q^{0}_{ij} \ ((i,j) \in F_{\ell}).  \label{eq:initialize31}
\end{eqnarray}
Obviously,  P$_{\ell}(\coneK)$ coincides 
with $\overline{\rm COP}(\coneK)$. Hence \eqref{eq:equivalence31} holds for $q = \ell$. 
As the induction hypothesis,  we assume 
that for $r \in \{\ell,\ldots, 2\}$, the objective coefficients $Q^{0q}_{ij} \in \Real$ $((i,j) \in F_q)$ 
of P$_q(\coneK)$ $(q=\ell,\ell-1,\ldots,r)$ have been computed so that \eqref{eq:equivalence31} holds for $q = \ell,\ldots,r$.
To show how $Q^{0(r-1)}_{ij} \in \Real$ $((i,j) \in F_{r-1})$ of P$_{r-1}(\coneK)$ are computed 
so that \eqref{eq:equivalence31} holds for $q = r-1$, 
we consider the following subproblem of P$_{r}(\coneK)$:
\begin{eqnarray}
\mbox{$\widetilde{\rm P}_{r}(\coneK^{C_r},\lambda)$:}
\quad \tilde{\eta}_r(\coneK^{C_r},\lambda) 
= \inf \left\{ G_r(\X_{C_r}) : 
 \X_{C_r} \in \Psi_r(\coneK^{C_r},\lambda) 
\right\}, 
\label{eq:Prr}
\end{eqnarray}
where $G_r(\X_{C_r}) :=$
$ \displaystyle 
\inprod{\Q^{0r}_{C_r}}{\X_{C_r}} -  Q^{0r}_{k_rk_r}X_{k_rk_r} 
= 
\displaystyle \sum_{(i,j) \in (C_r \times C_r) \backslash \{(k_r,k_r)\} } Q^{0r}_{ij}X_{ij}$,  
 and $\lambda \geq 0$ denotes a parameter. 
By~\eqref{eq:recursive},   we observe that 
\begin{eqnarray}
\eta_r(\coneK) & = & \inf \left\{ 
\displaystyle 
\sum_{(i,j) \in F_{r-1}} Q^{0r}_{ij}X_{ij} + G_r(\X_{C_r}) : 
\begin{array}{l}
(\X_{C_1},\ldots,\X_{C_{r-1}}) \in \Phi_{r-1}(\coneK), \\ \X_{C_r} \in \Psi_r(\coneK^{C_r},X_{k_rk_r}) 
\end{array}
\right\} \nonumber \\ 
& = & \inf \biggl\{  \ \displaystyle 
\sum_{(i,j) \in F_{r-1}} Q^{0r}_{ij}X_{ij} + \inf\left\{ 
 G_r(\X_{C_r}) : 
\X_{C_r} \in \Psi_r(\coneK^{C_r},X_{k_rk_r})  \right\}:  \nonumber \\
& &  \hspace{75mm} (\X_{C_1},\ldots,\X_{C_{r-1}}) \in \Phi_{r-1}(\coneK)\ \biggl\} \nonumber \\
& = & \inf \left\{  \ \displaystyle 
\sum_{(i,j) \in F_{r-1}} Q^{0r}_{ij}X_{ij} + 
\tilde{\eta}_r(\coneK^{C_r},X_{k_rk_r}) :  (\X_{C_1},\ldots,\X_{C_{r-1}}) \in \Phi_{r-1}(\coneK)\right\} 
\label{eq:subproblem31} \\ 
& = & \displaystyle \sum_{(i,j) \in F_{r-1}} Q^{0r}_{ij}\overline{X}_{ij} + 
\tilde{\eta}_r(\coneK^{C_r},\overline{X}_{k_rk_r}) \nonumber \\
& &  \hspace{40mm} \mbox{for every optimal solution $\overline{\X}$ of COP$(\coneK)$}. 
\label{eq:subproblem32} 
\end{eqnarray}
Here the last equality follows from the induction assumption. 

We now focus on the inner problem $\widetilde{\rm P}_{r}(\coneK^{C_r},X_{k_rk_r})$ 
with the optimal value $\tilde{\eta}_r(\coneK^{C_r},X_{k_rk_r})$.  Define
\begin{eqnarray}
\widetilde{\X}^{C_r} = \left\{
\begin{array}{ll}
\mbox{an optimal solution of $\widetilde{\rm P}_{r}(\coneK^{C_r},1)$}  & \mbox{if it exists}, \\ 
\O \in \SymMat^{C_r}  & \mbox{otherwise.}
\end{array} 
\right. \label{eq:widetildeX} 
\end{eqnarray}

\lemm  \label{lemma:correlative31} 
For a fixed $r$, 
assume that~\eqref{eq:equivalence31} holds for $q=\ell,\ldots,r$. 
Let  $(\X^*_{C_1},\ldots,\X^*_{C_{r-1}}) \in \Phi_{r-1}(\coneK)$ be an optimal solution of 
COP~\eqref{eq:subproblem31}, {\it i.e.}, 
there exists an optimal solution 
$\X_{C_r}^*$ of $\widetilde{\rm P}(\coneK^{C_r},X^*_{k_rk_r})$ with the optimal value 
$\tilde{\eta}(\coneK^{C_r},X^*_{k_rk_r})$ such that 
\begin{eqnarray*}
\displaystyle \sum_{(i,j) \in F_{r-1}} Q^{0r}_{ij}X^*_{ij} + \tilde{\eta}_r(\coneK^{C_r},X^*_{k_rk_r}) = \eta_r(\coneK). 
\end{eqnarray*}
Then $\X_{C_r} = X^*_{k_rk_r} \widetilde{\X}_{C_r}$ is 
an optimal solution of $\widetilde{\rm P}_{r}(\coneK^{C_r},X^*_{k_rk_r})$ with the 
optimal value $X^*_{k_rk_r} \tilde{\eta}_r(\coneK^{C_r},1)$,  
where $0 \times \tilde{\eta}_r(\coneK^{C_r},1) = 0$ is assumed 
even when $\tilde{\eta}_r(\coneK^{C_r},1) = \infty$, and $(\X^*_{C_1},\ldots,\X^*_{C_{r-1}},X^*_{k_rk_r} \widetilde{\X}_{C_r})$ is an optimal solution of P$_r(\coneK)$ with the optimal value $\eta_r(\coneK)$. 
\elemm 
\proof{
We consider two cases separately: $X^*_{k_rk_r} = 0$ and $X^*_{k_rk_r} > 0$. 
First, assume that $X^*_{k_rk_r} = 0$. In this case, $\X'_{C_r}$ is a feasible solution of 
$\widetilde{\rm P}_{r}(\coneK^{C_r},0)$ iff $\lambda\X'_{C_r}$ is a feasible solution of 
$\widetilde{\rm P}_{r}(\coneK^{C_r},0)$ for every $\lambda > 0$. This implies that $\tilde{\eta}_r(\coneK^{C_r},0)$ is 
either $0$ or $-\infty$. It follows from $\eta_r(\coneK) > -\infty$ that $\tilde{\eta}_r(\coneK^{C_r},0) = - \infty$ 
cannot occur. Hence $\tilde{\eta}_r(\coneK^{C_r},0) = 0 = X^*_{k_rk_r} \tilde{\eta}_r(\coneK^{C_r},1)$ follows 
(even when $\tilde{\eta}_r(\coneK^{C_r},1) = \infty$). 
We also see that $\X_{C_r} = \O =  X^*_{k_rk_r} \widetilde{\X}_{C_r}$ is a trivial optimal solution with the 
objective value $0$. 
Now assume that $X^*_{k_rk_r} > 0$. 
By assumption, $\widetilde{\rm P}_{r}(\coneK^{C_r},X^*_{k_rk_r})$ has an optimal solution $\X^*_{C_r}$. 
It is easy to see that 
$\X'_{C_r}$ is an optimal solution 
of $\widetilde{\rm P}_{r}(\coneK^{C_r},X^*_{k_rk_r})$ with the objective value 
$\xi$ iff $\X'_{C_r}/X^*_{k_rk_r}$ is an optimal solution of $\widetilde{\rm P}_{r}(\coneK^{C_r},1)$ with the objective value $\xi/X^*_{k_rk_r}$. 
We also see from the definition that $\widetilde{\X}_{C_r}$ is an optimal solution of 
 $\widetilde{\rm P}_{r}(\coneK^{C_r},1)$ with the objective value $\tilde{\eta}_r(\coneK^{C_r},1)$.  Hence 
$X^*_{k_rk_r} \widetilde{\X}_{C_r}$ is an optimal solution of  
$\widetilde{\rm P}_{r}(\coneK^{C_r},X^*_{k_rk_r})$ with the objective value 
$X^*_{k_rk_r} \tilde{\eta}_r(\coneK^{C_r},1)$.  Finally, we observe that 
$(\X^*_{C_1},\ldots,\X^*_{C_{r-1}},X^*_{k_rk_r} \widetilde{\X}^{C_r})$ is a feasible solution of P$_r(\coneK)$ with the objective value $\eta_r(\coneK)$. Therefore it is an optimal solution of P$_r(\coneK)$.}
\qed

By Lemma~\ref{lemma:correlative31}, $\tilde{\eta}_r(\coneK^{C_r},X_{k_rk_r})$ in  \eqref{eq:subproblem31} can be replaced 
with $X_{k_rk_r}\tilde{\eta}_r(\coneK^{C_r},1)$ and 
$\tilde{\eta}_r(\coneK^{C_r},\overline{X}_{k_rk_r})$  in \eqref{eq:subproblem32} 
with $\overline{X}_{k_rk_r}\tilde{\eta}_r(\coneK^{C_r},1)$.  Therefore, by defining  
\begin{eqnarray}
Q^{0r-1}_{ij} = \left\{
\begin{array}{ll}
Q^{0r}_{k_rk_r} + \tilde{\eta}_{r}(\coneK^{C_r},1) & \mbox{if } \ (i,j) = (k_r,k_r), \\[5pt] 
Q^{0r}_{ij} & \mbox{otherwise},  
\end{array} 
\right. \label{eq:update31}
\end{eqnarray}
we obtain 
\begin{eqnarray*}
\eta_r(\coneK) & = & \inf \left\{  \ \displaystyle 
\sum_{(i,j) \in F_{r-1}} Q^{0r-1}_{ij}X_{ij}  :  (\X_{C_1},\ldots,\X_{C_{r-1}}) \in \Phi_{r-1}(\coneK)\right\} \\
& = & \eta_{r-1}(\coneK) 
\ = \ \displaystyle \sum_{(i,j) \in F_{r-1}} Q^{0r-1}_{ij}\overline{X}_{ij} \ 
\mbox{for every optimal solution $\overline{\X}$ of COP$(\coneK)$}. 
\end{eqnarray*}
(We note that if $Q^{0r-1}_{kk} + \tilde{\eta}_r(\coneK^{C_r},1) = \infty$ occurs in~\eqref{eq:update31}, then  
 $\widetilde{\X}_{C_r}$ is set to be $\O \in \SymMat^{C_r}$). 
Thus we have shown that~\eqref{eq:equivalence31} holds for $q=r-1$ under the assumption that~\eqref{eq:equivalence31} holds for $q=\ell,\ldots,r$.  

Consequently, we obtain the following theorem by induction with decreasing $r$ from $\ell$ to $2$. 
\theo \label{theorem:main42}
Let $\coneK \in \{\bGamma^n,\CPP^n,\DNN^n\}$. 
Initialize the objective coefficients $Q^{0r}_{ij}$ $((i,j) \in F_r)$ of the sequence  P$_r(\coneK)$
$(r=\ell,\ldots,1)$ by~\eqref{eq:initialize31} for $r=\ell$, and update them by~\eqref{eq:update31} for $r=\ell,\ldots,2$. 
Then each P$_r(\coneK)$ in the sequence is equivalent to $\overline{\rm COP}(\coneK)$, more precisely, 
\eqref{eq:equivalence31} holds for $q=\ell,\ldots,1$. 
\etheo

\subsection{An algorithm for solving \mbox{$\overline{\rm COP}(\coneK)$}}

By Theorem~\ref{theorem:main42}, 
we know that $\overline{\rm COP}(\coneK)$ in \eqref{eq:COPbar} is equivalent to 
$P_1(\coneK)$, {\it i.e.},  $\eta_1(\coneK) = \bar{\zeta}(\coneK)$. 
Hence, the optimal value 
$\bar{\zeta}(\coneK)$ of $\overline{\rm COP}(\coneK)$ can be obtained by solving  
$P_1(\coneK)$. 

Lemma~\ref{lemma:correlative31}  suggests how to retrieve an optimal solution of $\overline{\rm COP}(\coneK)$. 
Let $\X^*_{C_1}$ be an optimal solution of 
P$_1(\coneK)$. 
For $r \in \{2,\ldots,\ell\}$, assume that 
an optimal solution  $(\X^*_{C_1},\ldots,\X^*_{C_{r-1}})$ of P$_{r-1}(\coneK)$ has been computed.
 Let 
$\X^*_{C_r} = \widetilde{\X}_{C_r} X^*_{k_rk_r}$, which is an optimal solution of 
$\widetilde{\rm P}_r(\coneK^{C_r},X^*_{k_rk_r})$ with the objective value $\tilde{\eta}(\coneK^{C_r},1) X^*_{k_rk_r}$. 
By Lemma~\ref{lemma:correlative31}, 
$(\X^*_{C_1},\ldots,\X^*_{C_r})$ is an optimal solution of P$_{r}(\coneK)$. 
We can continue this
procedure until an optimal solution of P$_{\ell}(\coneK)$ is obtained.  Note that P$_{\ell}(\coneK)$ is equivalent to $\overline{\rm COP}(\coneK)$.

\algo \label{algo:correlative} \mbox{ \ } 
\begin{description}
\item{Step 1: } (Initialization for computing $Q^{0r}_{ij} \ ((i,j) \in F_r, r=\ell,\ldots,1)$ and $\widetilde{\X}_{C_r}$ 
$(r=\ell,\ldots,2)$) 
Let $r = \ell$ and $Q^{0r}_{ij} = Q^0_{ij}$ $((i,j) \in F_{\ell}$). \vspace{-2mm} 
\item{Step 2: } If $r = 1$, go to Step 4. Otherwise, choose  $k_r \in C_r$ 
such that $(C_{1}\cup\cdots\cup C_{r-1})\cap C_r = \{k_r\}$. 
Solve $\widetilde{\rm P}_{r}(\coneK^{C_r},1)$. If $\widetilde{\rm P}_{r}(\coneK^{C_r},1)$ is infeasible, 
let $\widetilde{\X}_{C_r} = \O$ and $\tilde{\eta}_r(\coneK^{C_r},1) = \infty$. Otherwise, let 
$\widetilde{\X}_{C_r}$ be an optimal solution of $\widetilde{\rm P}_{r}(\coneK^{C_r},1)$ with the optimal value 
$\tilde{\eta}_r(\coneK^{C_r},1)$. Define $Q^{0r-1}_{ij}$ $((i,j) \in F_{r-1})$ by~\eqref{eq:update31}.  \vspace{-2mm} 
\item{Step 3: } Replace $r$ by $r-1$ and go to Step 2. \vspace{-2mm} 
\item{Step 4: } (Initialization for computing an optimal solution of $\overline{\rm COP}(\coneK)$) 
Solve P$_1(\coneK)$.  Let $\bar{\zeta}(\coneK) = \eta_1(\coneK)$ 
and $\X^*_{C_1}$ be an optimal solution of P$_1(\coneK)$. Let $r = 1$. \vspace{-2mm} 
\item{Step 5: } If $r=\ell$, then output the optimal value $\bar{\zeta}(\coneK)$ and an optimal solution 
$(\X^*_{C_1},\ldots,\X^*_{C_{\ell}})$ of $\overline{\rm COP}(\coneK)$. Otherwise go to Step 6. 
\vspace{-2mm} 
\item{Step 6: } Replace $r$ by $r+1$. 
Let $\X^*_{C_r} = X^*_{k_rk_r} \widetilde{\X}_{C_r}$. Go to Step 5.  \vspace{-2mm} 
\end{description}
\ealgo

\rema In Steps 1 through 3 of Algorithm~\ref{algo:correlative}, the problems $\widetilde{\rm P}_r(\coneK^{C_r},1)$ 
are solved 
sequentially from $\ell$ to $2$. This sequential order has been determined by the running intersection 
property~\eqref{eq:rip31}, which is induced from a clique tree of the block-clique graph $G(N,E)$. 
Recall that  $G(N,E)$ represents 
 the aggregated and correlative sparsity of COP$(\coneK)$.  
If $C_{r}$ other than $C_{\ell}$  is  also a leaf node of 
the clique-tree, then the problem $\widetilde{\rm P}_r(\coneK^{C_{r}},1)$ can be solved independently from 
the other problems and the objective coefficient $Q^{0r}_{k_rk_r}$ can be updated to $Q^{0r-1}_{k_rk_r}$ 
by~\eqref{eq:update31},  where it is  assumed that having chosen $C_1$ as the root node of the clique tree 
naturally determines all leaf nodes. 
In fact,  the problems associated with leaf nodes of the clique-tree 
can be solved in parallel. Furthermore, if we remove (some of) those nodes from the clique tree after solving 
their associated problems and updating the coefficients of the corresponding objective function  by~\eqref{eq:update31}, then 
new leaf nodes may appear. Then, the procedure of solving the problems associated with those new leaf nodes  
in parallel and  updating the objective coefficients 
can be repeatedly applied until we solve P$_1(\coneK)$. 
\erema

\subsection{Equivalence of $\overline{\rm COP}(\coneK_1)$ and $\overline{\rm COP}(\coneK_2)$ for 
$\coneK_1, \ \coneK_2  \in \{\bGamma^n,\CPP^n,$ $\DNN^n\}$}

Let $\coneK_1, \ \coneK_2 \in \{\bGamma^n,\CPP^n,\DNN^n\}$. For each of $s=\{1,2\}$, 
initialize the objective coefficients $Q^{0r}_{ij}$ $((i,j) \in F_r)$ of the sequence of P$_r(\coneK_s)$
$(r=\ell,\ldots,1)$ by~\eqref{eq:initialize31} for $r=\ell$, and update them by~\eqref{eq:update31} for $r=\ell,\ldots,2$. 
Then $\widetilde{\rm P}_r(\coneK_{1}^{C_{r}},1) $ and $\widetilde{\rm P}_{r}(\coneK_2^{C_{r}},1)$ 
share the same linear equality and inequality constraints, although the cone constraints $\X_{C_r} \in \coneK_1^{C_r}$ and 
$\X_{C_r} \in \coneK_2^{C_r}$ differ. 
Moreover, if $r = \ell$, 
the objective coefficients   of
$\widetilde{\rm P}_r(\coneK_{1}^{C_{r}},1) $ and $\widetilde{\rm P}_{r}(\coneK_2^{C_{r}},1)$ are the same as $\Q^{0}_{C_{\ell}}$.
Thus, it is natural to investigate the equivalence between $\widetilde{\rm P}_{\ell}(\coneK_{1}^{C_{\ell}},1)$ and $\widetilde{\rm P}_{\ell}(\coneK_2^{C_{\ell}},1)$ or whether 
$\tilde{\eta}_{\ell}(\coneK_1^{C_{\ell}}) = \tilde{\eta}_{\ell}(\coneK_2^{C_{\ell}}) $ holds. 
Let $r \in \{\ell,\ldots,2\}$. Assume that the objective coefficients of
$\widetilde{\rm P}_{q}(\coneK_{1}^{C_{q}},1)$ and $\widetilde{\rm P}_{q}(\coneK_2^{C_{q}},1)$ are the same, and  they  are 
equivalent, {\it i.e.},  $\tilde{\eta}_{q}(\coneK_1^{C_{q}}) = \tilde{\eta}_{q}(\coneK_2^{C_{q}}) $ holds 
$(q=\ell,\ldots,r)$. Then $\widetilde{\rm P}_{r-1}(\coneK_{1}^{C_{r-1}},1)$ and $\widetilde{\rm P}_{r-1}(\coneK_2^{C_{r-1}},1)$ share common objective coefficients and constraints except for the cone constraints $\X_{C_{r-1}} \in \coneK_1^{C_r-1}$ and 
$\X_{C_{r-1}} \in \coneK_2^{C_r-1}$. Thus, 
 the pair of 
$\widetilde{\rm P}_{r}(\coneK_{1}^{C_r},1)$  and $\widetilde{\rm P}_{r}(\coneK_2^{C_r},1)$ can be  compared  for their equivalence 
recursively  from $r=\ell$ to $r=2$. When all the pairs are equivalent, we can conclude by Theorem~\ref{theorem:main42} that $\eta_{1}(\coneK_1^{C_{1}}) = \bar{\zeta}(\coneK_1)$ and 
$\eta_{1}(\coneK_2^{C_{1}}) =\bar{\zeta}(\coneK_2)$. We also see that P$_1(\coneK_1)$ and P$_1(\coneK_1)$ share 
common objective coefficients. 
Consequently, the question on whether $\overline{\rm COP}(\coneK_1)$ and $\overline{\rm COP}(\coneK_2)$ are equivalent 
is reduced to the question on whether the pair of $\widetilde{\rm P}_{r}(\coneK_{1}^{C_{r}},1)$  and $\widetilde{\rm P}_{r}(\coneK_2^{C_r},1)$  is equivalent $(r=\ell,\ldots,2)$ and whether the pair of 
P$_1(\coneK_1)$ and P$_1(\coneK_2)$ is equivalent. 

For the convenience of the subsequent discussion, we introduce the sequence of the following COPs:
\begin{eqnarray*}
\widetilde{P}'_r(\coneK^{C_r}): \quad \tilde{\eta}'_r(\coneK^{C_r})  
& = & \inf \left\{ \inprod{\Q^{0r}_{C_r}}{\X} : \X \in \Psi_r(\coneK^{C_r},1) \right\} \\
& = & \inf \left\{ \inprod{\Q^{0r}_{C_r}}{\X} : 
\begin{array}{l}
\X \in \coneK^{C_r},  \ X_{k_rk_r} = 1, \\
\inprod{\Q^p_{C_r}}{\X_{C_r}} = 0 \ (p \in I_{\rm eq}^r),\\
\inprod{\Q^p_{C_r}}{\X_{C_r}} \leq 0 \ (p \in I_{\rm ineq}^r)
\end{array}
\right\}. 
\end{eqnarray*}
$(r = \ell,\ldots,1)$.  For each $r \in \{\ell,\ldots,2\}$, 
 $\widetilde{P}'_r(\coneK^{C_r})$ and $\widetilde{P}_r(\coneK^{C_r},1)$  
share a common feasible region $\Psi_r(\coneK^{C_r},1)$ and their objective values differ by a constant $Q^{0r}_{k_rk_r}$ 
for every $\X_{C_r} \in \Psi_r(\coneK^{C_r},1)$, which shows that they are essentially the same.  For 
$r=1$, $\widetilde{P}'_1(\coneK^{C_1})$ coincides with P$_1(\coneK)$. 
Thus,  
$\widetilde{\rm P}_r(\coneK^{C_r},1)$ $(r=\ell,\ldots,2)$ and 
P$_1(\coneK)$ can be dealt with 
as $\widetilde{P}'_r(\coneK^{C_r})$ 
$(r=\ell,\ldots.1)$. In particular, 
the question  on 
the equivalence of $\overline{\rm COP}(\coneK_1)$ and $\overline{\rm COP}(\coneK_2)$ can be stated as 
whether the pair of  $\widetilde{P}'_{r}(\coneK_{1}^{C_{r}})$ and $\widetilde{P}'_{r}(\coneK_2^{C_r})$
are equivalent for $r=\ell,\ldots,1$. 

Summarizing the discussions above, we obtain the following results.
\theo \label{theorem:main44}
Let $\coneK_s \in \{\bGamma^n,\CPP^n,\DNN^n\}$ $(s=1,2)$ and $\coneK_1 \not= \coneK_2$. For each $s = 1,2$, 
initialize the objective coefficients $Q^{0r}_{ij}$ $((i,j) \in F_r)$ of the sequence  P$_r(\coneK_s)$
$(r=\ell,\ldots,1)$ by~\eqref{eq:initialize31} for $r=\ell$, and update them by~\eqref{eq:update31} for $r=\ell,\ldots,2$. 
Assume that $\tilde{\eta}'_r(\coneK_1^{C_r}) = \tilde{\eta}'_r(\coneK_2^{C_r})$ $(r = \ell,\ldots,1)$. 
Then $\bar{\zeta}(\coneK_1)  =  \bar{\zeta}(\coneK_2) $. 
\etheo 

\subsection{Equivalence of  $\widetilde{P}'_r(\coneK_1)$ and  $\widetilde{P}'_r(\coneK_2)$ 
for $\coneK_1, \ \coneK_2  \in \{\bGamma^n,\CPP^n,$ $\DNN^n\}$
}

\theo \label{theorem:main45} Let $r \in {\ell,\ldots,1}$ be fixed. \vspace{-2mm} 
\begin{description}
\item{(i) } Assume that the aggregated sparsity of the data matrices $\Q^{0r}_{C_r}$ and $\Q^p_{C_r}$
$(p \in I_{\rm eq}^r\cup I_{\rm ineq}^r)$ 
is represented by a block-clique graph $G(C,E_C)$ with the maximal cliques of size at most $4$. 
Then $\tilde{\eta}'_r(\CPP^{C_r}) = \tilde{\eta}'_r(\DNN^{C_r})$. \vspace{-2mm} 
\item{(ii) } Assume that 
$ 
I_{\rm ineq}^r = \emptyset, \ 
\Psi_r(\bGamma^{C_r},1)\not=\emptyset, \
 \inprod{\Q^{0r}_{C_r}}{\X} \geq 0 \ \mbox{if } \X_{C_r} \in \Psi_r(\bGamma^{C_r},0) \ \mbox{and } 
\Q^p_{C_r} \in \ \SymMat^{C_r}_+ + \SymN^{C_r}
\ \mbox{(the dual of $\DNN^{C_r}$)}.
$ 
Then $\tilde{\eta}'_r(\bGamma^{C_r}) = \tilde{\eta}'_r(\CPP^{C_r})$. \vspace{-2mm} 
\item{(iii) } If the assumptions in (i) and (ii) above are satisfied, then $\tilde{\eta}'_r(\bGamma^{C_r}) = \tilde{\eta}'_r(\CPP^{C_r}) =\tilde{\eta}'_r(\DNN^{C_r})$. \vspace{-2mm} 
\item{(iv) } Let $\widetilde{C}_r=C_r\backslash\{k_rk_r\}$. Assume that $\Q^{0r}_{\widetilde{C}_r} \in \SymMat^{\widetilde{C}_r}_+$, 
$\Q^{p}_{C_r} \in \SymMat^{C_r}_+$ $(p \in I_{\rm eq}^r)$ and $\Q^{p}_{\widetilde{C}_r} \in \SymMat^{\widetilde{C}_r}_+$ 
$(p \in I_{\rm ineq}^r)$. Then $\tilde{\eta}'_r(\bGamma^{C_r}) =\tilde{\eta}'_r(\CPP^{C_r}) = \tilde{\eta}'_r(\DNN^{C_r})$. \vspace{-2mm} 
\item{(v) } Assume that the data matrices $\Q^{0r}_{C_r}$ and $\Q^p_{C_r}$
$(p \in I_{\rm eq}^r\cup I_{\rm ineq}^r)$ are all diagonal. Then 
$\tilde{\eta}'_r(\bGamma^{C_r}) =\tilde{\eta}'_r(\CPP^{C_r}) = \tilde{\eta}'_r(\DNN^{C_r})$. 
\vspace{-2mm} 
\end{description}
\etheo 
\proof{(i), (ii) and (iv) follow directly from (ii) of Theorem~\ref{theorem:main31},  Lemma~\ref{lemma:nonconvexCOP} and 
Lemma~\ref{lemma:convexCOP}, respectively. (i) and (ii) imply (iii). 
(v) follows from the discussion in Remark \ref{remark:LP} (b).
\qed

\bigskip
We note that the aggregated sparsity of the updated objective coefficient $\Q^{0r}_{C_r}$ 
is the same as the the original objective coefficient $\Q^{0}_{C_r}$ 
since only the diagonal elements 
$Q^0_{k_rk_r}$ $(r=\ell,\ldots,2)$ are updated by~\eqref{eq:update31}. Thus the assumption in (i) 
can be verified from the original data matrices $\Q^{0}_{C_r}$ 
and $\Q^p_{C_r}$ $( p \in  I_{\rm eq}^r \cup I_{\rm ineq}^r)$  
before solving $\overline{\rm COP}(\DNN^n)(\coneK)$ by Algorithm~\ref{algo:correlative}. 
On the other hand, the assumptions ``$ \inprod{\Q^{0r}_{C_r}}{\X} \geq 0 \ \mbox{if } \X_{C_r} \in \Psi_r(\bGamma^{C_r},0)$'' 
in (ii) and $\Q^{0r}_{\widetilde{C}_r} \in \SymMat^{\widetilde{C}_r}_+$ in (iv) depend on $\Q^{0r}_{C_r}$ $(r=\ell,\ldots,1)$ in general. 
If we know that 
$\Psi_r(\coneK^{C_r},0) = \{\O\}$, then ``$ \inprod{\Q^{0r}_{C_r}}{\X} \geq 0 \ \mbox{if } \X_{C_r} \in \Psi_r(\bGamma^{C_r},0)$'' 
obviously holds. In this case, the assumptions in (ii) can be verified from the original data matrices.  On the other hand, each  element $Q^{0r}_{ij}$ of the 
matrix $\Q^{0r}_{\widetilde{C}_{r}}$ is determined recursively by~\eqref{eq:update31} such that
\begin{eqnarray}
\left. 
\begin{array}{ll}
Q^{0\ell}_{ij} = Q^{0}_{ij},  \\ 
[3pt]
\displaystyle Q^{0q-1}_{ij} = \left\{
\begin{array}{ll}
Q^{0q}_{k_qk_q} + \tilde{\eta}_q(\coneK^{C_q},1) &  \mbox{if } (i,j) = (k_q,k_q), \\ 
[3pt]
Q^{0q}_{ij}  & \mbox{otherwise} 
\end{array}
\right. \ (q=\ell,\ldots,r). 
\end{array} \right\} \label{eq:update32}
\end{eqnarray}
As a result, only some of the diagonal elements of $\Q^{0r}_{\widetilde{C}_{r}}$ can differ 
from $\Q^{0}_{\widetilde{C}_{r}}$. 
For example, if $\Q^{0}_{\widetilde{C}_{r}}$ is positive semidefinite and 
$\Q^0_{C_q} \in \SymN^{C_q}$ for all $q=\ell, \ldots,r+1$, then $\Q^{0r}_{\widetilde{C}_{r}}$ is 
guaranteed to be  positive semidefinite,  since then 
$\tilde{\eta}_r(\coneK^{C_q}) \geq 0$ $(q=\ell, \ldots,r+1)$. In such a case, the assumptions in (iv) can be verified from the original data matrices. By~\eqref{eq:update32}, if 
the data matrices $\Q^{0}_{C_r}$ and $\Q^p_{C_r}$
$(p \in I_{\rm eq}^r\cup I_{\rm ineq}^r)$ are all diagonal, then the assumption of (v) is satisfied.

%
\section{Examples of QOPs}

We present
 two examples of QOPs  that can be reformulated as their DNN relaxations. 
 The problem in Section~5.1 is a nonconvex QOP with linear and 
complementarity constraints, and  the one in Section~5.2 is  a partially convex QOP with quadratic inequality constraints. 
They  are constructed as follows. First,  choose a block-clique graph $G(N,E)$ with the maximal cliques 
$C_q$ $(q=1,\ldots,\ell)$. For the first example in Section~5, we use the block-clique graph in  Figure 1 (a), 
 and  for the second example in~Section 5.2,  the one in Figure 1 (b). As the block-clique graph $G(N,E)$ 
induces a clique tree  (see Figure 2 (a) and  see  Figure 2 (b),
 respectively),  we renumber its $\ell$ maximal cliques 
so that they can satisfy the running intersection property~\eqref{eq:rip21}. 

Using the definitions of $\Phi_r(\coneK)$, 
$\Psi_r(\coneK,\lambda)$ $(r=\ell,\ldots,1)$ and~\eqref{eq:recursive}, we can rewrite $\overline{\rm COP}(\coneK)$ as 
\begin{eqnarray*}
\bar{\zeta}(\coneK) & = & \eta_{\ell}(\coneK) \ 
= \ \inf\left\{ \sum_{(i,j) \in F_{\ell}}Q^0_{ij}X_{ij} : (\X_{C_1},\ldots,\X_{C_{\ell}}) \in \Phi_{\ell}(\coneK) \right\} \\
& = & \inf\left\{ \sum_{(i,j) \in F_{\ell}}Q^0_{ij}X_{ij} :  
\X_{C_1} \in \Phi_{1}(\coneK), \
\X_{C_r} \in \Psi(\coneK^{C_r},X_{k_rk_r}) \ (r = \ell,\ldots,2)
\right\} \\ 
& = & \inf\left\{ \sum_{(i,j) \in F_{\ell}}Q^0_{ij}X_{ij} : \X_{C_r} \in \Psi_r(\coneK^{C_r},X_{k_rk_r}) \ (r = \ell,\ldots,1), \ X_{11} = 1 \right\}
\end{eqnarray*}
Instead of $\overline{\rm COP}(\coneK)$ itself, we describe the problem with its subproblems 
\begin{eqnarray*}
\widehat{\rm P}_r(\coneK^{C_r},\lambda): \ \hat{\eta}_r(\coneK^{C_r},\lambda) = \inf \left\{\inprod{\Q^0_{C_r}}{\X_{C_r}}:  \X_{C_r} \in \Psi_r(\coneK^{C_r},\lambda)\right\}
\end{eqnarray*}
using 
\begin{eqnarray*}
C_r & : & k_r \in C_r, \ \Q^0_{C_r} \in \SymMat^{C_r} \ \mbox{(and its property if any)} \ \mbox{and } \Psi_r(\bGamma^{C_r},\lambda) 
\end{eqnarray*}
for $(r=1,\ldots,\ell)$, where $\lambda \geq 0$ denotes a parameter. 
Notice that $\widehat{\rm P}(\coneK^{C_r},\lambda)$ is similar to  
the problem $\widetilde{\rm P}_r(\coneK^{C_r},\lambda)$ introduced in Section~4.1.
We also note that  the objective coefficient $\Q^{0r}_{C_r}$ 
is updated from $\Q^{0}_{C_r}$ in the construction of the sequence P$_r(\coneK)$ $(r=\ell,\ldots,1)$ by~\eqref{eq:update31}, 
while the objective coefficient of $\widehat{\rm P}(\coneK^{C_r},\lambda)$ is fixed to $\Q^0_{C_r}$.
Indeed, the description $\widehat{\rm P}(\coneK^{C_r},\lambda)$ is sufficient and more convenient to 
execute Algorithm~\ref{algo:correlative} for 
computing the optimal value $\bar{\zeta}(\coneK)$ of $\overline{\rm COP}(\coneK)$ and also to apply Theorems~\ref{theorem:main44} and~\ref{theorem:main45} to establish 
$\zeta(\bGamma^n) = \zeta(\CPP^n) = \zeta(\DNN^n)$. 

\subsection{A QOP with linear and complementarity constraints}

Consider the block-clique graph $G(N,E)$ given in Figure 1 (a). In this case, $n=11$ and 
$N=\{1,\ldots,11\}$. 
To represent a QOP  as $\overline{\rm COP}(\bGamma^n)$, let 
\begin{eqnarray*}
C_1=\{1,2,3,4\} &:& k_1=1, \ 
\Q^0_{C_1} \in \SymMat^{C_1}, \\
& &  \Psi_1(\bGamma^{C_1},\lambda) = \left\{ \x_{C_1}\x_{C_1}^T : 
\begin{array}{l}
\x_{C_1} \in \Real^{C_1}_+, \ x_1= \lambda, \\
-8x_1+x_2+x_3+2x_4=0 
\end{array}
\right\}, \\
C_2=\{3,5,6,7\} &:& k_2=3, \ 
\Q^0_{C_2} \in \SymMat^{C_2}, \\
& &  \Psi_1(\bGamma^{C_2},\lambda) = \left\{ \x_{C_2}\x_{C_2}^T : 
\begin{array}{l}
\x_{C_2} \in \Real^{C_2}_+, \ x_3= \lambda, \  x_5x_6=0, \\
-x_3+x_5+2x_6+x_7=0
\end{array}
\right\}, 
\end{eqnarray*}
\begin{eqnarray*}
C_3=\{3,8,9\} &:& k_3=3, \ 
\Q^0_{C_3} \in \SymMat^{C_3}, \\
& &  \Psi_1(\bGamma^{C_3},\lambda) = \left\{ \x_{C_3}\x_{C_3}^T : 
\begin{array}{l}
\x_{C_3} \in \Real^{C_3}_+, \ x_3= \lambda, \\
-2x_3+3x_8+x_9=0 
\end{array}
\right\}, \\
C_4=\{5,10,11\} &:& k_4=5, \ 
\Q^0_{C_4} \in \SymMat^{C_4}, \\
& &  \Psi_1(\bGamma^{C_4},\lambda) = \left\{ \x_{C_4}\x_{C_4}^T : 
\begin{array}{l}
\x_{C_4} \in \Real^{C_4}_+, \ x_5= \lambda, \\
-x_5+x_{10}+x_{11}=0 
\end{array}
\right\}.
\end{eqnarray*}

The assumption in (i) of Theorem~\ref{theorem:main45} is satisfied since ${C_r}$ $(r=1,2,3,4)$ 
 are of size at most $4$. 
All assumptions in (ii) of Theorem~\ref{theorem:main45} are also satisfied. In fact, 
$I_{\rm ineq}^r = \emptyset$, $\Psi_r(\bGamma^{C_r},1) \not=\emptyset$ and $\Psi_r(\bGamma^{C_r},0) = \{\O\}$, which 
implies that  $\inprod{\Q^{0r}_{C_r}}{\X} \geq 0 \ \mbox{if } \X_{C_r} \in \Psi_r(\bGamma^{C_r},0)$, 
$(r=1,2,3,4)$. 
As shown in Example~\ref{example:bbc},  each linear equality constraint can be rewritten  
in $\Psi_r(\bGamma^{C_r},\lambda)$ 
$(r=1,2,3,4)$ as $\inprod{\Q_{C_r}}{\x_{C_r}\x_{C_r}^T} = 0$ for some $\Q_{C_r} \in \SymMat^{C_r}_+$ and 
the complementarity constraint $x_5x_6=0$ in $\Psi_2(\bGamma^{C_r},\lambda)$ as $\inprod{\Q_{C_r}}{\x_{C_r}\x_{C_r}^T} = 0$ 
for some $\Q_{C_r} \in \SymN^{C_r}$. By (i) and (ii) of Theorem~\ref{theorem:main45}, 
$\tilde{\eta}'_r(\bGamma^{C_r}) = \tilde{\eta}'_r(\CPP^{C_r}) =  \tilde{\eta}'_r(\DNN^{C_r})$ $(r=1,2,3,4)$. Therefore, 
by Theorem~\ref{theorem:main44}, it follows that $\bar{\zeta}(\bGamma^n) = \bar{\zeta}(\CPP^n) 
= \bar{\zeta}(\DNN^n)$ holds. 

\subsection{A partially convex QOP}

Consider the block-clique graph $G(N,E)$ given in  Figure 1 (b). In this case, $n=13$ and 
$N=\{1,\ldots,13\}$. To represent a QOP  as $\overline{\rm COP}(\bGamma^n)$, let 
\begin{eqnarray*}
C_1 = \{1,2,3\} & : &  k_1= 1, \ \Q^0_{C_1} \in \SymMat^{C_1}, \\
& &  \Psi_1(\bGamma^{C_1},\lambda) = \left\{ \x_{C_1}\x_{C_1}^T : 
\begin{array}{l}
\x_{C_1} \in \Real^{C_1}_+, \ 
x_1= \lambda, \\
-x_1+x_2+x_3 =0 
\end{array}
\right\}, \\
C_2 = \{3,4,5,6,7\} & : &  k_2= 3, \ \Q^p_{C_2} \in \SymMat^{C_2} : \mbox{diagonal } 
 \ (p = 0,1,2),\\
& &  \Psi_2(\bGamma^{C_2},\lambda) = \left\{ \x_{C_2}\x_{C_2}^T : 
\begin{array}{l}
\x_{C_2} \in \Real^{C_2}_+, \ x_3= \lambda,\\
  \inprod{\Q^1_{C_2}}{\x_{C_2}\x_{C_2}^T } = 0\\
  \inprod{\Q^2_{C_2}}{\x_{C_2}\x_{C_2}^T } \leq 0
\end{array}
\right\}  \not= \emptyset,  \\ 
C_3 = \{3,8,9,10,11\} & : &  k_3= 3, \
\Q^p_{C_3} \in \SymMat^{C_3}, \ \widetilde{C}_3 = \{8,9,10,11\}, \ 
 \Q^p_{\widetilde{C}_3} \in \SymMat^{\widetilde{C}_3}_+ \ (p = 0,3,4), \\
 & & \Psi_3(\bGamma^{C_3},\lambda) = \left\{ \x_{C_3}\x_{C_3}^T : 
\begin{array}{l}
\x_{C_3} \in \Real^{C_3}_+, \ x_3= \lambda, \\ -x_3-2x_8+x_9+x_{10} = 0, \\
 \inprod{\Q^p_{C_3}}{\x_{C_3}\x_{C_3}^T } \leq 0 \ (p=3,4)
\end{array}
\right\} \not= \emptyset, \\ 
C_4 = \{11,12,13\} & : &  k_4= 11, \ \Q^0_{C_4} \in \SymMat^{C_4}, \\
& &  \Psi_4(\bGamma^{C_4},\lambda) = \left\{ \x_{C_4}\x_{C_4}^T : 
\begin{array}{l}
\x_{C_4} \in \Real^{C_4}_+, \ x_{11} = \lambda, \\
 -x_{11}+x_{12}+2x_{13} = 0 
\end{array} 
\right\}.
\end{eqnarray*}

For $r=1,4$, we similarly see that  all assumptions in (i) and (ii) of 
Theorem~\ref{theorem:main45} are satisfied as in   the example in Section 5.1.
Thus, the identity 
$\tilde{\eta}_r(\bGamma^{C_r}) = \tilde{\eta}_r(\CPP^{C_r}) = 
\tilde{\eta}_r(\DNN^{C_r})$ holds $(r=1,4)$. 

Since the data matrices $\Q^p_{C_2} \in \SymMat^{C_2}$ $(p =0,1,2)$ are diagonal, 
$\widehat{\rm P}_2(\coneK^{C_2},\lambda)$ 
becomes an LP of the form 
\begin{eqnarray*}
\hat{\eta}_r(\coneK^{C_2},\lambda) = \inf \left\{\sum_{i\in C_2}Q^0_{ii}X_{ii}:  
\begin{array}{l}
X_{ii} \in \Real_+ \ (i \in C_2), \ X_{33} = \lambda, \\
\sum_{i\in C_2}Q^1_{ii}X_{ii} = 0, \ 
\sum_{i\in C_2}Q^1_{ii}X_{ii} \leq 0
\end{array}
\right\}. 
\end{eqnarray*}
Hence (v) of Theorem~\ref{theorem:main45}, $\tilde{\eta}_2(\bGamma^{C_2}) = \tilde{\eta}_2(\CPP^{C_2}) = \tilde{\eta}_2(\DNN^{C_2})$ holds. 

To see whether the identity  holds for $r=3$, we need to apply (iv) of Theorem~\ref{theorem:main45} 
since $\Psi_3(\coneK^{C_3},\lambda)$ involves convex quadratic inequality. 
It suffices to check whether $\Q^{03}_{\widetilde{C}_3}$ 
is positive semidefinite. 
By the assumption and the updating formula~\eqref{eq:update32}, we know that 
\begin{eqnarray*}
Q^{03}_{ij} = \left\{
\begin{array}{ll}
Q^{0}_{11,11} + \tilde{\eta}_4(\bGamma^{C_4}) & \mbox{if } (i,j) = (11,11), \\
[3pt]
Q^0_{ij} & \mbox{otherwise}. 
\end{array}
\right.
\end{eqnarray*}
For example, if $\tilde{\eta}_4(\bGamma^{C_4})$ is 
nonnegative, then 
$\Q^{03}_{\widetilde{C}_3} \in \SymMat^{\widetilde{C}_3}_+$ and the identity 
$\tilde{\eta}_3(\bGamma^{C_3}) = \tilde{\eta}_3(\CPP^{C_3})= 
\tilde{\eta}_3(\DNN^{C_3})$ holds. 
Therefore, 
by Theorem~\ref{theorem:main44}, it follows that $\bar{\zeta}(\bGamma^n) = \bar{\zeta}(\CPP^n) 
= \bar{\zeta}(\DNN^n)$ holds. 

%
\section{Concluding remarks}

Nonconvex QOPs and CPP problems are known to be NP hard and/or numerically intractable in general,
as opposed to computationally tractable
DNN problems.
Thus, finding some classes of QOPs 
or CPP problems that are equivalent to DNN problems  is an essential problem
in the study of the  theory and applications of nonconvex QOPs.
Two major obstacles to finding such classes are: (A) $\CPP^n$ is a proper subset of $\DNN^n$ if 
$n \geq 5$ and (B) general quadratically constrained nonconvex QOPs are NP hard and  
numerically intractable. 
As a result, CPP reformulations of a class of QOPs with linear equality, binary and complementarity 
constraints in nonnegative variables   still remain numerically intractable.  
One way to overcome these obstacles is to ``decompose" the cone $\CPP^n$ 
 into cones with size at most 4 and/or to ``decompose'' a QOP into 
convex QOPs of any size and linearly constrained nonconvex QOPs with variables at most 4. To obtain such decompositions,  
a fundamental method is exploiting  structured sparsity. 

To describe a QOP, its CPP and DNN relaxations, 
COP$(\coneK)$ with $\coneK \in \{\bGamma^n,\CPP^n,$ $\DNN^n\}$ has been introduced in Section 1. 
In Section 3, we have provided a method to decompose 
the cone $\CPP^n$ of the CPP relaxation, COP$(\CPP^n)$ of the QOP described as COP$(\bGamma^n)$
into cones with size at most 4 by exploiting the aggregated sparsity of the data 
matrices $\Q^p$ $(p \in \{\0\}\cup I_{rm eq} \cup I_{\rm ineq})$ of COP$(\CPP^n)$ which have been 
 represented by a block-clique graph $G(N,E)$. In Section 4, a method to decompose the QOP itself into 
convex QOPs of any size and linearly constrained nonconvex QOPs with variables at most 4 has been presented  by exploiting their correlative sparsity. 

As for the structured sparsity that leads to the equivalence among COP$(\bGamma^n)$, COP$(\CPP^n)$ and 
COP$(\DNN^n)$, the aggregated and/or correlative sparsity of the data matrices
represented with a block-clique graph have played a crucial role in our discussion. 
We should mention that block-clique graphs may not be frequently observed in a wide class of QOPs.
It is interesting, however, to construct a new optimization model 
based on the structure provided by a block-clique graph. 

For further development of such an optimization model, we emphasize that if a QOP described as COP$(\bGamma^n)$ can be solved exactly then it can be incorporated in the model as a subproblem.
In \cite{KIM2003}, it was shown that the exact solutions of
nonconvex QOPs with nonnpositive off-diagonal data matrices  can be found. The result has been applied to the
optimal power flow problems \cite{LOW2014a}.  
More precisely,
assume that $I_{\rm eq} = \emptyset$ and that all off-diagonal elements of 
$\Q^p$ $(p \in \{0\} \cup I_{\rm ineq})$ are nonpositive. Then the QOP described as COP$(\bGamma^n)$ can be solved exactly by its 
SDP and SOCP relaxation \cite[Theorem 3.1]{KIM2003}. Thus, the QOP can be incorporated in our model as a subproblem. 
The details are omitted here.


\end{document}